\documentclass[11pt]{article}

\usepackage{amsmath, amssymb, amscd}
\usepackage{epsfig}

\def\eqref#1{(\ref{#1})}
\newcommand{\goth}{\frak}
\newcommand{\g}{{\frak g}}
\newcommand{\arrow}{{\:\longrightarrow\:}}
\newcommand{\Z}{{\Bbb Z}}
\newcommand{\C}{{\Bbb C}}
\newcommand{\R}{{\Bbb R}}

\newcommand{\6}{\partial}
\newcommand{\1}{\sqrt{-1}\:}
\newcommand{\restrict}[1]{{\left|_{{\phantom{|}\!\!}_{#1}}\right.}}

\renewcommand{\c}[1]{{\cal #1}}
\newcommand{\calo}{{\cal O}}

\renewcommand{\tilde}{\widetilde}
\renewcommand{\bar}{\overline}
\renewcommand{\phi}{\varphi}
\renewcommand{\epsilon}{\varepsilon}
\renewcommand{\geq}{\geqslant}
\renewcommand{\leq}{\leqslant}

\newcommand{\End}{\operatorname{End}}

\newcommand{\Id}{\operatorname{Id}}

\newcommand{\Hom}{\operatorname{Hom}}
\newcommand{\Ext}{\operatorname{Ext}}

\newcommand{\codim}{\operatorname{codim}}

\newcommand{\coker}{\operatorname{coker}}

\newcommand{\rk}{\operatorname{rk}}

\newcommand{\Def}{\operatorname{Def}}
\newcommand{\Tw}{\operatorname{Tw}}
\newcommand{\Tr}{\operatorname{Tr}}
\newcommand{\Spec}{\operatorname{Spec}}

\newcommand{\comment}[1]{{}}

\def\blacksquare{\hbox{\vrule width 4pt height 4pt depth 0pt}}
\def\endproof{\blacksquare}

\makeatletter

\@ifundefined{Bbb}
     {\newcommand{\Bbb}[1]{{\mathbb #1}}}%
{}

\newcommand{\ps@verbit}{%
  \renewcommand{\@oddhead}{%
          \scriptsize
          {Coherent sheaves on generic K3 surfaces and tori}
          \hfil\tiny {M. Verbitsky, \ \ \ \ \ May 17, 2002 }}
  \renewcommand{\@evenhead}{\@oddhead}
  \renewcommand{\@oddfoot}{\hfil\thepage\hfil}
  \renewcommand{\@evenfoot}{\@oddfoot}}
 
\newcounter{Mycounter}[section]
\newcounter{lemma}[section]
\renewcommand{\thelemma}{{Lemma \thesection.\arabic{lemma}}}
\newcommand{\lemma}{%
     \setcounter{lemma}{\value{Mycounter}}
     \refstepcounter{lemma}
     \stepcounter{Mycounter}
     {\bf \thelemma:\ }}

\newcounter{claim}[section]
\renewcommand{\theclaim}{{Claim \thesection.\arabic{claim}}}
\newcommand{\claim}{%
     \setcounter{claim}{\value{Mycounter}}
     \refstepcounter{claim}
     \stepcounter{Mycounter}
     {\bf \theclaim:\ }}

\newcounter{sublemma}[section]
\newcounter{corollary}[section]
\renewcommand{\thecorollary}{{Corollary \thesection.\arabic{corollary}}}
\newcommand{\corollary}{%
     \setcounter{corollary}{\value{Mycounter}}
     \refstepcounter{corollary}
     \stepcounter{Mycounter}
     {\bf \thecorollary:\ }}

\newcounter{theorem}[section]
\renewcommand{\thetheorem}{{Theorem \thesection.\arabic{theorem}}}
\newcommand{\theorem}{%
     \setcounter{theorem}{\value{Mycounter}}
     \refstepcounter{theorem}
     \stepcounter{Mycounter}
     {\bf \thetheorem:\ }}

\newcounter{conjecture}[section]
\newcounter{proposition}[section]
\renewcommand{\theproposition}
       {{Proposition \thesection.\arabic{proposition}}}
\newcommand{\proposition}{%
     \setcounter{proposition}{\value{Mycounter}}
     \refstepcounter{proposition}
     \stepcounter{Mycounter}
     {\bf \theproposition:\ }}

\newcounter{definition}[section]
\renewcommand{\thedefinition}
       {{Definition~\thesection.\arabic{definition}}}
\newcommand{\definition}{%
     \setcounter{definition}{\value{Mycounter}}
     \refstepcounter{definition}
     \stepcounter{Mycounter}
     {\bf \thedefinition:\ }}

\newcounter{example}[section]
\newcounter{remark}[section]
\renewcommand{\theremark}{{Remark \thesection.\arabic{remark}}}
\newcommand{\remark}{%
     \setcounter{remark}{\value{Mycounter}}
     \refstepcounter{remark}
     \stepcounter{Mycounter}
     {\bf \theremark:\ }}

\newcounter{problem}[section]
\newcounter{question}[section]
\@addtoreset{equation}{section}
\@addtoreset{footnote}{section}
\makeatother

\begin{document}

\begin{center}
{\LARGE\bf
Coherent sheaves on general K3 surfaces and tori
}
\\[4mm]
Misha Verbitsky,\footnote{The author is 
partially supported by CRDF grant RM1-2354-MO-02.}
\\[4mm]
{\tt verbit@thelema.dnttm.ru, \ \  verbit@mccme.ru}
\end{center}


{\small 
\hspace{0.15\linewidth}
\begin{minipage}[t]{0.7\linewidth}
{\bf Abstract} \\
Let $M$ be a K3 surface
or an even-dimensional  compact torus. We show
that the category of coherent sheaves
on $M$ is independent from the choice of the
complex structure, if this complex structure is generic.

\end{minipage}
}

{
\small
\tableofcontents
}


\section{Introduction}
\label{_Intro_Section_}


\subsection{An overview}

A hyperk\"ahler manifold 
is a  Riemannian manifold $M$ equipped with an orthogonal
action of quaternion algebra in its tangent space $TM$,
in such a way that for all $L\in {\Bbb H}$, $L^2=-1$,
the $L$ induces a K\"ahler structure on $M$.

Let $M$ be a hyperk\"ahler manifold, and $I\in \Bbb H$
a quaternion such that $I^2=-1$. By $(M,I)$ we understand
the $M$ considered as a K\"ahler manifold, with the complex structure
$I$. 

We study the category of coherent sheaves on $(M,I)$.
This question can be stated more algebraically as follows. 
Using the linear algebra of the quaternionic action, it is easy
to construct a holomorphic symplectic form on $(M,I)$,
that is, a holomorphic, closed, nowhere degenerate $(2,0)$-form.
Conversely, every compact holomorphic symplectic K\"ahler
manifold admits, by Calabi-Yau theorem, 
a unique hyperk\"ahler metrics in a given 
K\"ahler class. Details of this correspondence
can be found e.g. in \cite{_Besse:Einst_Manifo_}.

The subject of this paper can be understood algebraically as
the study of coherent sheaves on compact holomorphically 
symplectic K\"ahler manifolds. 

A holomorphic symplectic manifold has trivial canonical class;
the canonical class is trivialized by taking the top exterior
power of the holomorphic symplectic form. Mirror conjecture,
in the form proposed by M. Kontsevich (\cite{_Kontsevich:ICM_}), 
relates the category of coherent sheaves on a Calabi-Yau manifold
and the Gromov-Witten geometry of the Lagrangian cycles on 
its mirror partner. More precisely, this version of the
mirror conjecture states that the derived category
of coherent sheaves should be equivalent to the
Fukaya category of the mirror dual manifold.

The mirror partner of a holomorphic symplectic manifold
is again holomorphic symplectic. However, a generic deformation of 
a holomorphic symplectic manifold does not admit holomorphic curves,
so one would expect that all Gromov-Witten-type invariants
(including the Fukaya category) are (in some sense)
trivial. 

In this paper we prove the mirror counterpart of this conjecture,
for $M$ a K3 surface or a compact complex torus. 

\hfill

\theorem \label{_main_intro_Theorem_}
Let $M_1$, $M_2$ be K3 surfaces or 
compact complex tori of dimension $2d$. Assume that
$M_1$, $M_2$ are generic
in the sense of having no non-trivial integer 
$(p,p)$-cycles, for $0<p<2d$, and $\c C_1$
$\c C_2$ be the category of coherent sheaves on $M_1$, $M_2$.
Then $\c C_1$ is (non-canonically) equivalent to $\c C_2$.

\hfill

{\bf Proof:} See \ref{_generic_are_equi_Theorem_}.
\endproof

\hfill

In the algebraic situation, 
D. Orlov (\cite{_Orlov:K3_})  
studied how the derived 
category of coherent sheaves 
behaves when one changes 
the complex structure 
of a K3 surface. This is not the situation
we work in, because generic complex structures
on hyperk\"ahler manifolds are never algebraic.

The derived category of coherent sheaves 
is a weaker invariant than the 
category  of coherent sheaves itself.
Using the Fourier-Mukai transform,
D. Orlov shows that the derived
categories of coherent sheaves on the projective 
K3 surfaces $M$ and $M'$ are equivalent
if and only if the Hodge lattices of trancendental
cycles in $H^2(M)$, $H^2(M')$ are
isometric.

\subsection{A plan of the proof}

The proof of \ref{_main_intro_Theorem_} is based on hyperk\"ahler
geometry and on the theory of reflexive sheaves. A coherent sheaf
is called {\bf reflexive} if it is isomorphic to its second dual
(\ref{_refle_Definition_}).

Let $M$ be a compact hyperk\"ahler manifold, and $S$ the
2-sphere of complex structures induced by the quaternionic action.
There is a natural $SU(2)$-action on the cohomology $H^*(M)$
of $M$ (\ref{_SU(2)_commu_Laplace_Lemma_}). 
A complex structure $L\in S$ is called
{\bf generic} (\ref{_generic_manifolds_Definition_}) 
if all $(p,p)$-cycles on $(M,L)$ are $SU(2)$-invariant. 

In \cite{_Verbitsky:Symplectic_II_} it was proven that 
all induced complex structures $L\in S$ outside of a countable
set are generic (\ref{_generic_are_dense_Proposition_}). 

If $L$ is generic, then all integer $(1,1)$-classes
are $SU(2)$-invariant, hence are of degree zero
(\ref{_Lambda_of_inva_forms_zero_Lemma_}). 
Therefore, all vector bundles
on $(M,L)$ are semistable of degree 
zero.\footnote{Throughout 
this paper, stability, polystability 
and semistability of vector bundles and coherent sheaves
is understood in the sense of Mumford-Takemoto, see 
\ref{_degree,slope_destabilising_Definition_}.
The degree of a coherent sheaf $F$ on 
a K\"ahler manifold $M$ is equal to
$\int_M c_1(F)\wedge \omega^{n-1}$, where $n=\dim_\C M$ and
$\omega$ is the K\"ahler form.}

\hfill

\theorem\label{_bundle_equi_Theorem_}
Let $M$ be a hyperk\"ahler K3 surface or 
a compact hyperk\"ahler 
torus of real dimension $2d$.  Let
$L_1$, $L_2$ induced complex structures, and 
$\c C_1^r$, $\c C_2^r$ the categories of 
reflexive sheaves on $(M,L_1)$
and $(M, L_2)$ (\ref{_refle_category_Definition_}). Assume that $(M, L_1)$ 
and $(M, L_2)$ have no non-trivial integer 
$(p,p)$-cycles, for $0<p <d$.\footnote{This is
equivalent to $(M,L_i)$ being Mumford-Tate generic;
see Section \ref{_conne_hype_Section_} for details.}
Then the categories $\c C_1^r$, $\c C_2^r$ 
are equivalent. 

\hfill

{\bf Proof:} This is \ref{_refle_shea_isomo_Theorem_}. \endproof

\hfill

The proof of \ref{_bundle_equi_Theorem_} is based 
on the following idea, which is due to
\cite{_NHYM_}. 
Let $L$ be a generic induced complex structure on $M$.
Every bundle on $(M, L)$ is semistable.
In \cite{_NHYM_} we studied the twistor correspondence,
which allows one to construct a Yang-Mills
bundle $\Tw(B)$ on the twistor space from a 
stable bundle $B$ on $(M, L)$ (see \ref{_twi_functor_Theorem_}). 
This correspondence is actually injective: 
restricting $\Tw(B)$ to
$(M,L)\subset \Tw(M)$, one
obtains $B$ again. In \cite{_NHYM_} 
we outline the set of conditions which 
guarantees that a given holomorphic bundle
on the twistor space is obtained from
a bundle with connection on $(M,L)$ 
(see \ref{_twi_functor_Theorem_}). 
This formalism is functorial; we call it
the direct and inverse twistor transform.

In this paper we apply the twistor formalism to
the semistable bundles. Given a semistable bundle
on $(M, L)$, we lift it to a holomorphic bundle
on a twistor space, thus obtaining an equivalence between
the category of semistable bundles on $(M, L)$ and a certain
subcategory of holomorphic bundles on the twistor space.

\hfill

If we replace the category $\c C^{b}$
of holomorphic vector bundles
by its subcategory $\c C^{st}$
of polystable bundles\footnote{A 
polystable bundle (\ref{_degree,slope_destabilising_Definition_}) 
is by definition a direct sum of stable bundles of the same slope.},
then \ref{_bundle_equi_Theorem_}  is well known
(see \cite{_V:Hyperholo_sheaves_}). The
equivalence $\c C_1^{st}\cong \c C_2^{st}$ is 
canonical.

\hfill

The sphere $S^2$ of induced complex structures gives a holomorphic
embedding $\C P^1 \stackrel j \hookrightarrow \c M$, where $\c M$ is the
moduli of complex structures on $M$. Let 
$I_1$, $I_2\in j(\C P^1)$ be generic
induced complex structures. By \ref{_bundle_equi_Theorem_}, 
the corresponding categories of holomorphic vector bundles are
isomorphic. In \cite{_coho_announce_} it was proven that
any pair of complex K\"ahler structures $I_1$, $I_2\in \c M$
within the same component of the moduli space can be connected
by a sequence of rational curves
of form $\C P^1 \stackrel {j_k} \hookrightarrow \c M$
associated with a sequence of hyperk\"ahler structures.
Moreover, these hyperk\"ahler structures can be chosen 
in such a way that the sequential intersection points
$j_k(\C P^1)\cap j_{k+1}(\C P^1)$ correspond to generic
complex structures (\ref{_generic_are_connected_Theorem_}).

Therefore, \ref{_bundle_equi_Theorem_} shows that the category
of reflexive sheaves on $(M,I)$ is independent from 
the complex structure $I$, as long as $I$ remains generic.

\hfill

\subsection{Flat bundles on complex tori}

Let $T$ be a generic compact 
hyperk\"ahler torus of complex 
dimension $>2$, and $I$ a generic induced
complex structure. Then, 
all coherent sheaves on $(T,I)$  have isolated
singularities and locally free reflexizations
(\ref{_cohe_on_toru_Theorem_}).

Denote by $\c C^b(T,I)$ the category of
holomorphic vector bundles
on $(T, I)$. To prove that $\c C^b(T,I)$
is independent from the choice of generic induced
complex structure $I$, we produce a (non-Hermitian)
flat connection $\nabla_B$ on a given bundle $B$
on $(T, I)$, in a functorial way. 
For $B$ stable, $\nabla$ coinsides with the 
Yang-Mills connection produced by the 
Uhlenbeck-Yau theorem (\ref{_UY_Theorem_}).

Let $I'$ be an arbitrary induced complex structure.
Taking the $(0,1)$-part of $\nabla_B$ with
respect to $I'$, we obtain a holomorphic 
bundle $B_{I'}$ over $(T, I')$.
We show that this gives an equivalence
between the respective categories
of holomorphic vector bundles $\c C^b(M,I)$ 
and $\c C^b(M,I')$, if $I'$ is also generic. 

The singularities are dealt with in the same fashion.
For any $x\in M$, we identify the neighbourhood
$U_I$, $U_{I'}$ of $x$ in $(M, I)$ and $(M, I')$
(Subsection \ref{_twi_spa_local_Subsection_}). 
We show that this identification is compatible
with the connection $\nabla_B$, in such a way
that the restriction of $B$ to $U_I$
is equivalent to the restriction
of $B_{I'}$ to $U_{I'}$
(Subsection \ref{_twi_spa_local_Subsection_}). This
allows to ``glue in'' the isolated singularities
to $B$ and $B_{I'}$ in a way which is
canonical and compatible with $\nabla_B$.

\hfill

The twistor correspondence (\ref{_twi_functor_Theorem_}) 
gives a way to speak of connections on $B$ algebraically:
we identify the connections with $SU(2)$-\-inva\-riant curvature
and the lifts of the holomorphic structure from $B$
to the twistor space. Throughout this paper, we speak
of holomorphic structures on the bundles
in the twistor space, but \ref{_twi_functor_Theorem_}
claims that this is the same as to speak of connections.

If we work over a torus, flatness is a more natural
condition than $SU(2)$-invariance of the curvature. 
We propose a way to prove \ref{_bundle_equi_Theorem_}
on a torus without using the hyperk\"ahler structure,
as follows.

\hfill

Let $(B,\nabla)$ be a flat bundle on a complex manifold.
Taking the $(0,1)$-part of $\nabla$, we obtain a
holomorphic structure operator on $B$. This
gives a functor $\tau$ from the category of
flat bundles to the category of holomorphic bundles.
In this paper, we construct the ``inverse'' functor,
mapping holomorphic bundles to flat ones.

\hfill

\claim\label{_holo_to_flat_torus_Claim_}
Let $T$ be a compact $d$-dimensional complex torus, 
without non-trivial integer $(p,p)$-cycles in
cohomology, $0<p<d$, and $\c C^b(T)$ the category of holomorphic
vector bundles on $T$. Assume that $d>2$. Then there exists
a functor $B \arrow (B, \nabla_B)$ from 
$\c C^b(T)$ to the category of flat bundles on $T$,
such that for any $B$, the holomorphic bundle $\tau(B, \nabla)$
is isomorphic to $B$, and for $B$ stable, 
the connection $\nabla$ is Hermitian. 

\hfill

{\bf Proof:}
This is \ref{_bun_admits_flat_Corollary_}. \endproof

\hfill

It seems that it is possible 
to deduce the equivalence of the categories
of vector bundles from
\ref{_holo_to_flat_torus_Claim_}. 
We take a vector bundle $B$ and 
obtain the flat connection $\nabla$ on $B$ as in
\ref{_holo_to_flat_torus_Claim_}.
We change the complex structure $I$ on $T$
to a sufficiently close $I'$ and restrict 
$(B, \nabla)$ to $(T, I')$,
obtaining a holomorphic bundle $B_{I'}$
on $(T,I')$. Since the stable bundles
correspond to flat Hermitian connections,
this construction will identify the
stable bundles on $(T,I)$ and $(T,I')$.
To prove that we obtained an 
equivalence, it remains to show that
this functor preserves the $\Ext^1$-groups
(see \ref{_embe_on_Ext_1_hence_full_faith_Corollary_}). 
However, the $\Ext$ classes can be represented by 
parallel forms on $(T, I)$, $(T,I')$,
and as we change the complex structure
from $I$ to $I'$, parallel forms remain 
parallel, and the corresponding cohomology classes
are non-zero, if $I'$ is sufficiently
close to $I'$.

To prove the equivalence of the respective 
categories of coherent sheaves, we use 
\ref{_cohe_on_toru_Theorem_}, stating that on
a generic compact complex torus, all 
coherent sheaves have isolated singularities
and smooth reflexive hulls. To pass
from a bundle to a coherent sheaf, we 
need to ``glue in'' isolated singularities
to a holomorphic vector bundle $B$.

Using the flat structure obtained in
\ref{_holo_to_flat_torus_Claim_},
we can introduce flat coordinates on 
the tori $(T,I)$, $(T,I')$ in a neighbourhood
of a given point $x$, in such a way that the
flat structure on $B$ will be compatible with these coordinates.
This gives a natural functorial isomorphism
between $B$ and $B_{I'}$ in a neighbourhood of $x$. 
As in the proof of \ref{_tilde_C(I)_C_I_equiv_Theorem_},
we could ``glue in'' the singularities to $B$ and $B_{I'}$
simultaneously, obtaining the isomorphism of 
the categories of coherent sheaves.

\subsection{Reflexive sheaves on hyperk\"ahler manifolds}

The main motivation of \ref{_bundle_equi_Theorem_}
is the following. Suppose that $X$ is a complex manifold
without divisors. We identify 
the category of reflexive sheaves 
$\c C^r(X)$ with the category
\[ \lim_\arrow \c C^b(X\backslash S), \]
where the limit is taken over all increasing
sequences of closed complex analytic subvarieties
$S\subset X$, $\codim_X S \geq 2$, 
and $\c C^b(X\backslash S)$ is the category
of holomorphic vector bundles on $X \backslash S$
which can be extended to $X$ as coherent sheaves
(\ref{_refle_category_Definition_}). 
From this description it is clear that
 every object of $\c C^r(X)$ is a 
finite extension of a sequence of simple objects.
We say that an abelian category with this property
has {\bf finite length}, or {\bf satisfies the
ascending and the descending chain conditions}.

Now consider a compact hyperk\"ahler manifold $M$ and let $I$ be a 
generic complex structure induced by quaternions.
All complex subvarieties of $(M,I)$ are hyperk\"ahler,
hence have even codimension (\cite{_Verbitsky:Symplectic_II_}). 
Therefore, $(M, I)$ has no divisors. The set of simple
objects of $\c C^r(M,I)$ is independent from $I$
as long as $I$ stays generic (\cite{_V:Hyperholo_sheaves_}).
The $\Ext^*$ groups between simple objects can also be computed
(see \ref{_direct_ima_twi_Proposition_}),
and (at least for locally free sheaves)
these groups are independent from $I$ as well.
We arrive at the following situation, which
was known in many instances since 1990-ies.
Given a hyperk\"ahler manifold $M$ ,
and generic complex structures $I_1$, $I_2$ in the same
connected component of the moduli space, 
consider the categories of reflexive sheaves
$\c C^r(M,I_1)$ and $\c C^r(M,I_2)$.
Then 
\begin{description}
\item[(i)]
These categories are of finite length, 
that is, satisfy the ascending and the descending chain conditions.
\item[(ii)] There is
a natural equivalence $R$ between the simple objects of
$\c C^r(M,I_1)$ and $\c C^r(M,I_2)$ (see
\cite{_V:Hyperholo_sheaves_}).
\item[(iii)] If $M$ is a K3 surface or a torus, then
the differential graded algebra
$\Ext^*(B,B)$ is isomorphic to $\Ext^*(R(B),R(B))$,
for any simple object $B$ in $\c C^r(M,I_1)$.
\end{description}
The statement (iii) is known only for vector bundles
(\cite{_Verbitsky:Hyperholo_bundles_}). Fortunately,
reflexive sheaves on $M$ 
are locally free, if $M$ is a generic torus
or a K3 surface (see \ref{_cohe_on_toru_Theorem_}). 

It seems that the conditions (i)-(iii) are already sufficient to prove
the equivalence of the categories $\c C^r(M,I_1)$ and $\c C^r(M,I_2)$,
so one could prove \ref{_bundle_equi_Theorem_}
in an algebraic fashion via the theory of abelian categories.
In the present paper we work in a different direction. 
We construct
a subcategory $\c C^{tw}_I$ of the category $\c C^r(\Tw(M))$
of reflexive sheaves on the twistor space of $M$
(see Subsection \ref{_C^tw_I_defini_Subsection_}).
We show that the natural restriction functor
from $\Tw(M)$ to $(M,I_i)$ induces
an isomorphism between $\c C^{tw}_I$ and $\c C^r(M,I_i)$.
This implies that the categories 
$\c C^r(M,I_1)$ and $\c C^r(M,I_2)$ are equivalent.
Varying the hyperk\"ahler structure, we obtain
the same result for all generic 
complex structures.

\subsection{Contents}

This paper is organized as follows.

\begin{itemize}

\item This Introduction explains some ideas 
used in this paper, in a heuristic fashion.

\item In Sections \ref{_hk_Section_} - \ref{_Hyperholomo_Section_}
we recall some results and definitions 
from the literature. We relate some basic
facts from hyperk\"ahler geometry, Yang-Mills
theory and geometry of stable bundles.

\item In Section \ref{_conne_hype_Section_},
we study the moduli of complex structures on
hyperk\"ahler manifolds. We define the Mumfort-Tate 
generic complex structures and show that these can
be connected by a sequence of hyperk\"ahler structures.

\item In Section \ref{_refle_Section_}
we recall some standard results from the theory
of reflexive sheaves (see e.g. \cite{_OSS_}).

\item In Section \ref{_holo_and_twi_Section_} we deal
with the twistor formalizm introduced in \cite{_NHYM_}.
Let $M$ be a compact hyperk\"ahler manifold, and $I$
a generic induced complex structure. 
Given a stable bundle $B$ on $(M,I)$,
we construct a holomorphic vector bundle $\Tw(B)$ on its twistor
space. We compute the cohomology of $\Tw(B)$ it terms of $H^*(B)$.

\item In Section \ref{_twi_local_Section_}
we study the local geometry of the twistor space $\Tw(M)$,
using the rational curves on $\Tw(M)$. Given a bundle $F$ on $\Tw(M)$
obtained from twistor transform, we trivialize $F$ 
in a neighbourhood of $(\C P^1 \backslash I) \subset \Tw(M)$,
where $\C P^1 = s_x$ is the curve of form
\[ \C P^1 \times \{x\} \subset \C P^1 \times M = \Tw(M)
\]
and  $I\subset s_x$ is an arbitrary point.
We study the associated trivializations on 
the sheaves of differential operators and local cohomology.
In particular, we prove that a differential operator or a 
section of the local cohomology sheaf  in
a neighbourhood of $s_x \backslash I$ can be extended
to the whole $s_x$ if it is compatible with the trivialization. 

\item In Section \ref{_tori_Section_}
we study the coherent sheaves on a  generic 
compact complex torus $T$, $\dim_\C \geq 3$. We show that
 all reflexive sheaves on $T$ are locally free.
We construct a natural flat connection
on any holomorphic bundle on $T$. 

\item In Section \ref{_Yoneda_Section_}
we deal with formal properties of abelian categories
satisfying the ascending and descending chain conditions (such 
categories are also called {\bf of finite length}).
We define the abelian categories of cohomological dimension $\leq 1$
and give a set of criteria which imply that a given
category has cohomological dimension $\leq 1$.
Given a functor of  categories of finite length
and cohomological dimension $\leq 1$, we show that
it is an equivalence of categories 
if it is an equivalence on the 
simple objects and on $\Ext^1$-groups.

\item In Section \ref{_refle_equi_Section_} we show that the categories
of reflexive sheaves on $(M, I)$ and $(M, I')$ are 
equivalent if $M$ is a generic K3 surface or a hyperk\"ahler
torus, and $I$, $I'$ are 
generic induced complex structures. 
For a surface, we use the results about 
local cohomology obtained earlier. For $M$
a torus of $\dim_\C >2$, we use the flat connection
constructed in Section \ref{_tori_Section_}
on every holomorphic bundle.

\item In Section \ref{_cohe_isola_Section_},
we apply the local geometry of the twistor 
space to ``glue in'' the singularities to 
vector bundles in a canonical way. This is used to show
that the category of coherent sheaves with isolated singularities
and smooth reflexizations on $(M, I)$ is independent
from the choice of induced complex structure $I$,
if $I$ remains generic.

\end{itemize}


\section{Hyperk\"ahler manifolds}
\label{_hk_Section_}


In this section we reproduce 
well known results from hyperk\"ahler geometry, for later use.
We follow \cite{_Besse:Einst_Manifo_} and 
\cite{_Verbitsky:Symplectic_II_}.

\subsection{Hyperk\"ahler manifolds and quaternionic action}

\definition \label{_hyperkahler_manifold_Definition_} 
(\cite{_Besse:Einst_Manifo_}) A {\bf hyperk\"ahler manifold} is a
Riemannian manifold $M$ endowed with three complex structures $I$, $J$
and $K$, such that the following holds.
 
\begin{description}
\item[(i)]  the metric on $M$ is K\"ahler with respect to these complex 
structures and
 
\item[(ii)] $I$, $J$ and $K$, considered as  endomorphisms
of a real tangent bundle, satisfy the relation 
$I\circ J=-J\circ I = K$.
\end{description}

\hfill 

The notion of a hyperk\"ahler manifold was 
introduced by E. Calabi (\cite{_Calabi_}).

\hfill

Clearly, a hyperk\"ahler manifold has a natural action of
the quaternion algebra ${\Bbb H}$ in its real tangent bundle $TM$. 
Therefore its complex dimension is even.
For each quaternion $L\in \Bbb H$, $L^2=-1$,
the corresponding automorphism of $TM$ is an almost complex
structure. It is easy to check that this almost 
complex structure is integrable (\cite{_Besse:Einst_Manifo_}).

\hfill

\definition \label{_indu_comple_str_Definition_} 
Let $M$ be a hyperk\"ahler manifold, and $L$ a quaternion satisfying
$L^2=-1$. The corresponding complex structure 
on $M$ is called
{\bf an induced complex structure}. The $M$, considered as a K\"ahler
manifold, is denoted by $(M, L)$. In this case,
the hyperk\"ahler 
structure is called {\bf compatible
with the complex structure $L$}.

\hfill

Let $M$ be a hyperk\"ahler manifold. We identify the group $SU(2)$
with the group of unitary quaternions. This gives a canonical 
action of $SU(2)$ on the tangent bundle, and all its tensor
powers. In particular, we obtain a natural action of $SU(2)$
on the bundle of differential forms. 

\hfill

The following lemma is clear.

\hfill

\lemma \label{_SU(2)_commu_Laplace_Lemma_}
The action of $SU(2)$ on differential forms commutes
with the Laplacian.
 
{\bf Proof:} This is Proposition 1.1
of \cite{_Verbitsky:Symplectic_II_}. \endproof
 
\hfill

Thus, for compact $M$, we may speak of the natural action of
$SU(2)$ in cohomology.

\hfill

Further in this article, we use the following statement.

\hfill

\lemma \label{_SU(2)_inva_type_p,p_Lemma_} 
Let $\omega$ be a differential form over
a hyperk\"ahler manifold $M$. The form $\eta$ is $SU(2)$-invariant
if and only if it is of Hodge type $(p,p)$ with respect to all 
induced complex structures on $M$.

\hfill

{\bf Proof:} Let $I$ be an induced complex structure,
and ${\goth u}_I:\; U(1) \arrow SU(2)$ the corresponding
embedding, induced by the map 
$\R= \goth{u}(1) \arrow \goth{su}(2)$, $1\arrow I$.
The Hodge decomposition on $\Lambda^*(M)$ 
coinsides with the weight decomposition of 
the $U(1)$-action ${\goth u}_I$.
An $SU(2)$-invariant form is
also invariant with respect to ${\goth u}_I$,
and therefore has Hodge type $(p,p)$. Conversely,
if a $\eta$ is invariant with respect to ${\goth u}_I$,
for all induced complex structures $I$, then
$\eta$ is invariant with respect to the
Lie group $G$ generated by these 
$U(1)$-subgroups of $SU(2)$.
A trivial linear-algebraic 
argument ensures that the group $G$ 
coinsides with the whole $SU(2)$.
This proves \ref{_SU(2)_inva_type_p,p_Lemma_}.
\endproof

\subsection{Generic induced complex structures 
and trianalytic subvarieties}

Let $M$ be a compact hyperk\"ahler manifold, $\dim_\R M =2m$.
 
\hfill

\definition\label{_trianalytic_Definition_} 
Let $N\subset M$ be a closed subset of $M$. Then $N$ is
called {\bf trianalytic} if $N$ is a complex analytic subset 
of $(M,L)$ for any induced complex structure $L$.
 
\hfill
 
Let $I$ be an induced complex structure on $M$,
and $N\subset(M,I)$ be
a closed analytic subvariety of $(M,I)$, $dim_\C N= n$.
Consider the homology class 
represented by $N$. Let $[N]\in H^{2m-2n}(M)$ denote 
the Poincare dual cohomology class, so-called
{\bf fundamental class} of $N$. Recall that
the hyperk\"ahler structure induces the action of 
the group $SU(2)$ on the space $H^{2m-2n}(M)$.
 
\hfill

\theorem\label{_G_M_invariant_implies_trianalytic_Theorem_} 
Assume that $[N]\in  H^{2m-2n}(M)$ is invariant with respect
to the action of $SU(2)$ on $H^{2m-2n}(M)$. Then $N$ is 
 trianalytic.
 
{\bf Proof:} This is Theorem 4.1 of 
\cite{_Verbitsky:Symplectic_II_}.
\endproof

\hfill
 
\remark \label{_triana_dim_div_4_Remark_}
Trianalytic subvarieties have an action of the quaternion algebra in
the tangent bundle. In particular,
the real dimension of such subvarieties is divisible by 4.

\hfill

\definition \label{_generic_manifolds_Definition_} 
Let $M$ be a compact hyperk\"ahler manifold,
and $I$ an induced complex
structure. We say that $I$ is {\bf of general type}
or {\bf generic} with respect to 
the hyperk\"ahler structure on $M$,
if all elements of the group
\begin{equation}\label{_gene_Equation_} 
  \bigoplus\limits_p H^{p,p}(M)\cap H^{2p}(M,\Z)\subset H^*(M)
\end{equation}
are $SU(2)$-invariant.

\hfill

\proposition \label{_generic_are_dense_Proposition_} 
Let $M$ be a compact hyperk\"ahler manifold,  and $S$
the set of induced complex structures over $M$. Denote by 
$S_0\subset S$ the set of generic induced complex structures.
Then $S_0$ is dense in $S$. Moreover, the complement
$S\backslash S_0$ is countable.

{\bf Proof:} This is Proposition 2.2 from
\cite{_Verbitsky:Symplectic_II_}
\endproof

\hfill

\ref{_G_M_invariant_implies_trianalytic_Theorem_} has the following
immediate corollary:

\hfill

\corollary \label{_hyperkae_embeddings_Corollary_} 
Let $M$ be a compact hyperk\"ahler manifold,
and $I$ a generic induced complex
structure.
Let $X\subset (M,I)$ be a closed complex analytic
subvariety. Then $X$ is trianalytic. 
In particular, $X$ has even complex dimension

\endproof

\hfill

Throughout this paper, except Section
\ref{_tori_Section_}, a weaker form of 
``generic'' can be used: we could say that $I$ is generic
if all integer $(p,p)$-cycles
are $SU(2)$-invariant, for $p=1,2$.


\section[Hyperholomorphic bundles]{Hyperholomorphic bundles}
\label{_Hyperholomo_Section_}


\subsection{Hyperholomorphic connections }

This Subsection contains several versions of a
definition of hyperholomorphic connection in a complex
vector bundle over a hyperk\"ahler manifold.
We follow \cite{_Verbitsky:Hyperholo_bundles_}.

 Let $B$ be a holomorphic vector bundle over a complex
manifold $X$, $\nabla$ a  connection 
in $B$ and $\Theta\in\Lambda^2\otimes End(B)$ be its curvature. 
This connection
is called {\bf compatible with the holomorphic structure} if
$\nabla_\gamma(\zeta)=0$ for any holomorphic section $\zeta$ and
any antiholomorphic tangent vector field $\gamma\in T^{0,1}(X)$. 
If there exists a holomorphic structure compatible with the given
Hermitian connection then this connection is called
{\bf integrable}.

\hfill
 
One can define the {\bf Hodge decomposition} in the space of differential
forms with coefficients in any complex bundle, in particular,
$End(B)$.

\hfill

\theorem \label{_Newle_Nie_for_bu_Theorem_}
Let $\nabla$ be a Hermitian connection in a complex vector
bundle $B$ over a complex manifold $X$. Then $\nabla$ is integrable
if and only if $\Theta\in\Lambda^{1,1}(X, \End(B))$, where
$\Lambda^{1,1}(X, \End(B))$ denotes the forms of Hodge
type (1,1). Also,
the holomorphic structure compatible with $\nabla$ is unique.

{\bf Proof:} This is Proposition 4.17 of \cite{_Kobayashi_}, 
Chapter I.
\endproof

\hfill

This proposition is a version of Newlander-Nirenberg theorem.
For vector bundles, it was proven by M. Atiyah and R. Bott.

\hfill

\definition \label{_hyperho_conne_Definition_}
Let $B$ be a Hermitian vector bundle with
a connection $\nabla$ over a hyperk\"ahler manifold
$M$. Then $\nabla$ is called {\bf hyperholomorphic} if 
$\nabla$ is
integrable with respect to each of the complex structures induced
by the hyperk\"ahler structure. 
 
As follows from 
\ref{_Newle_Nie_for_bu_Theorem_}, 
$\nabla$ is hyperholomorphic
if and only if its curvature $\Theta$ is of Hodge type (1,1) with
respect to any of the complex structures induced by a hyperk\"ahler 
structure.

As follows from \ref{_SU(2)_inva_type_p,p_Lemma_}, 
$\nabla$ is hyperholomorphic
if and only if $\Theta$ is an $SU(2)$-invariant differential form.

\subsection{Hyperholomorphic bundles, Yang-Mills connections and stability}

\definition\label{_degree,slope_destabilising_Definition_} 
Let $F$ be a coherent sheaf over
an $n$-dimensional compact K\"ahler manifold $M$. We define
{\bf the degree} $\deg(F)$ (sometimes the degree
is also denoted by $\deg c_1(F)$) as
\[ 
   \deg(F)=\int_M\frac{ c_1(F)\wedge\omega^{n-1}}{vol(M)}
\] 
and $\text{slope}(F)$ as
\[ 
   \text{slope}(F)=\frac{1}{\text{rk}(F)}\cdot \deg(F). 
\]
The number $\text{slope}(F)$ depends only on a
cohomology class of $c_1(F)$. 

Let $F$ be a coherent sheaf on $M$
and $F'\subset F$ its subsheaf 
with $0<\rk F'<\rk F$. Then $F'$ is 
called {\bf destabilizing subsheaf} 
if $\text{slope}(F') \geq \text{slope}(F)$

A coherent sheaf $F$ is called {\bf  stable}
\footnote{In the sense of Mumford-Takemoto},
or $\mu$-stable, 
if it has no destabilizing subsheaves. 
A coherent sheaf $F$ is called {\bf 
polystable} if it is a direct sum 
of stable sheaves of the same slope.
A coherent sheaf $F$ is called {\bf 
semistable} if for all
destabilizing subsheaves $F'\subset F$
we have $\text{slope}(F') = \text{slope}(F)$

\hfill

Let $M$ be a K\"ahler manifold with a K\"ahler form $\omega$.
Consider the standard Hodge operator 
on differential forms, $L: \; \eta\arrow\omega\wedge\eta$.
There is also a fiberwise-adjoint Hodge operator $\Lambda$
(see \cite{_Griffi_Harri_}).

\hfill

\definition \label{_Yang-Mills_Definition_}
Let $B$ be a holomorphic bundle over a K\"ahler manifold $M$
with a holomorphic Hermitian connection $\nabla$ and a 
curvature $\Theta\in\Lambda^{1,1}\otimes End(B)$.
The Hermitian metric on $B$ and the connection $\nabla$
defined by this metric are called {\bf 
Yang-Mills} if 

\begin{equation} \label{_YM_Equation_}
   \Lambda(\Theta)=c \cdot \Id\restrict{B},
\end{equation}
where $\Lambda$ is a Hodge operator, $c$ a constant, 
and $\Id\restrict{B}$ is 
the identity endomorphism which is a section of $End(B)$.

\hfill

Clearly, the constant $c$ is proportional to the
slope of $B$. Throughout this paper, we shall consider
only bundles of slope zero. In this case, the Yang-Mille
equation can be written simly as $\Lambda(\Theta)=0$.

\hfill

The following fundamental 
theorem provides examples of 
Yang-\-Mills \linebreak bundles.

\hfill

\theorem \label{_UY_Theorem_} 
(Uhlenbeck-Yau)
Let B be a
holomorphic bundle over a compact K\"ahler manifold. Then $B$ admits
a Hermitian Yang-Mills connection if and only if it is 
polystable. Moreover, the Yang-Mills 
connection is unique, if it exists.
 
{\bf Proof:} \cite{_Uhle_Yau_}. \endproof

\hfill

\remark\label{_harmonic_curv_Remark_}
It is easy to see that a connection is Yang-Mills
if and only if its curvature is harmonic.

\hfill

\proposition \label{_hyperholo_Yang--Mills_Proposition_}
Let $M$ be a hyperk\"ahler manifold, $L$
an induced complex structure and $B$ be a complex vector
bundle over $(M,L)$. 
Then every 
hyperholomorphic connection $\nabla$ in $B$
is Yang-Mills and satisfies $\Lambda(\Theta)=0$,
where $\Theta$ is a curvature of $\nabla$.
 
\hfill

{\bf Proof:} We use the definition of a hyperholomorphic 
connection as one with $SU(2)$-invariant curvature. 
Then \ref{_hyperholo_Yang--Mills_Proposition_}
follows from the following elementary observation

\hfill

\lemma \label{_Lambda_of_inva_forms_zero_Lemma_}
Let $\Theta\in \Lambda^2(M)$ be a $SU(2)$-invariant 
differential 2-form on $M$. Then
$\Lambda_L(\Theta)=0$ for each induced complex structure
$L$.\footnote{By $\Lambda_L$ we understand the Hodge operator 
$\Lambda$ associated with the K\"ahler complex structure $L$.}

{\bf Proof:} This is Lemma 2.1 of \cite{_Verbitsky:Hyperholo_bundles_}.
\endproof
 
\hfill

Let $M$ be a compact hyperk\"ahler manifold, 
and $I$ an induced 
complex structure. For any 
stable holomorphic bundle on $(M, I)$ there exists a unique
Hermitian  Yang-Mills connection 
which, for some bundles, turns out to be hyperholomorphic. 
It is possible to tell exactly when
this happens.

\hfill

\theorem \label{_inva_then_hyperho_Theorem_}
Let $B$ be a 
polystable holomorphic bundle over
$(M,I)$, where $M$ is a hyperk\"ahler manifold and $I$
is an induced complex structure over $M$. Then 
$B$ admits a  
hyperholomorphic connection if and only
if it is polystable and 
the first two Chern classes $c_1(B)$ and $c_2(B)$ are 
$SU(2)$-invariant.\footnote{We use \ref{_SU(2)_commu_Laplace_Lemma_}
to speak of action of $SU(2)$ in cohomology of $M$.}

{\bf Proof:} This is Theorem 2.5 of
 \cite{_Verbitsky:Hyperholo_bundles_}. \endproof

\hfill

\definition\label{_hh_polysta_Definition_}
Let $M$ be a compact hyperk\"ahler manifold,
$I$ an induced complex structure, and $B$ a polystable
bundle on $(M, I)$. We say that $B$ is hyperholomorphic
if $B$ admits a hyperholomorphic connection; 
equivalently, $B$ is hyperholomorphic
if the first two Chern classes $c_1(B)$ and $c_2(B)$ are 
$SU(2)$-invariant.

\hfill

{}From the definition of generic induced complex structures,
we immediately obtain the following corollary

\hfill

\corollary \label{_gene_hh_Corollary_}
Let $M$ be a compact hyperk\"ahler manifold,
$I$ an induced complex structure, which 
is generic in the sense of \ref{_generic_manifolds_Definition_},
and $B$ a stable bundle over $(M, I)$. Then 
$B$ is hyperholomorphic. 

\endproof


\section[Generic complex structures on hyperk\"ahler manifolds]{Generic complex structures \\on hyperk\"ahler manifolds}
\label{_conne_hype_Section_}


Let $M$ be a compact hyperk\"ahler manifold.
In this Section, we follow the arguments of \cite{_coho_announce_}
(see also \cite{_V:Hyperholo_sheaves_}), which were used
to establish the equivalence between various 
algebro-geometric structures on $(M, I)$
and $(M, J)$, where $I$ and $J$ are generic
complex structures on $M$ induced by
a hyperk\"ahler structure in the
same connected component of the moduli space.

Given a complex structure $I$ on $M$,
let $\c I:\; H^i(M) \arrow H^i(M)$
be the operator mapping a $(p, q)$-class
$\eta$ to $\1(p-q)\eta$. This gives a $U(1)$-action 
${\goth u}_I$ on $H^i(M)$. 

\hfill

\remark\label{_type_p_p_u_i_inv_Remark_}
Clearly, a form is
of type $(p, p)$ with respect to $I$ 
if and only if it is  ${\goth u}_I$-invariant. 

\hfill

Let $\g_0\subset \End(H^*(M))$ 
be the Lie algebra generated by ${\goth u}_I$,
for all complex structures $I$ which are compatible
with some hyperk\"ahler structure on $M$. 
In \cite{_coho_announce_}, this Lie algebra was computed
explicitly, for $M$ with 
$h^1(M)=0$, $h^{2,0}(M)=1$ (such 
hyperk\"ahler manifolds are called {\bf simple}). 
We have shown that $\g_0\cong \goth{so}(H^2(M), h)$,
where $h$ is a certain non-degenerate symmetric form on $H^2(M)$,
called {\bf the Bogomolov-Beauville form}. 

\hfill

By  \ref{_type_p_p_u_i_inv_Remark_},
a cohomology class $\eta$ is $\g_0$-invariant
if and only if $\eta$ is invariant with respect 
to all complex structures $I$ induced by some
hyperk\"ahler structure. For any 
cohomology class $\eta$ which 
is not $\g_0$-invariant, let $C_\eta$ be the
set of all complex structures $I$ which
satisfy $\goth u_I(\eta) =0$. Clearly,
$C_\eta$ is a closed complex analytic
subset in the moduli $\c M$ 
of complex structures on $M$.

Let $S\subset H^*(M, \Z)$ be the set
of all cohomology classes which 
are not $\g_0$-invariant,
and $C_S$ the union of $C_\eta$ for
all $\eta \in S$. By definition,
$C_S$ is a countable union of complex subvarieties
of positive codimension. Therefore, it has measure
zero, and its complement is dense in the moduli 
space $\c M$. We obtain the following definition
of generic (\cite{_V:Hyperholo_sheaves_}).

\hfill

\definition\label{_M-T_generic_Definition_}
Let $M$ be a compact manifold of hyperk\"ahler
type, and $I$ a complex structure which is
compatible with some hyperk\"ahler structure.
Consider the Lie algebra $\g_0\subset \End(H^*(M))$  defined above.
Assume that all integer $(p, p)$-classes on $(M, I)$ are 
$\g_0$-invariant. Then $I$ is called Mumford-Tate generic.

\hfill

\remark
When $M$ is a K3 surface or a torus, $M$ is Mumford-Tate generic
if and only if $M$ does not have non-trivial integer $(p,p)$-cycles.

\hfill

The following claim is clear from the definition.

\hfill

\claim\label{_MT_gene_gene_wrt_hk_Claim_}
Let $M$ be a hyperk\"ahler manifold, and $I$ an induced
complex structure. Assume that $I$ is Mumford-Tate
generic. Then $I$ is generic, in the sense of 
\ref{_generic_manifolds_Definition_}.

\endproof

\hfill

The following theorem is proven in \cite{_coho_announce_}
(see also \cite{_V:Hyperholo_sheaves_}).

\hfill

\theorem\label{_generic_are_connected_Theorem_}
Let $M$ be a compact hyperk\"ahler manifold,
and $I$, $I'$ complex structures of K\"ahler type 
which belong to the same connected component 
of moduli space. Then there exists a sequence 
$\c H_1, \c H_2, ... \c H_n$ of
hyperk\"ahler structures, and a set
$I_0 = I, I_2, I_3, ..., I_{n+1} = I'$
of complex structures on $M$, such that
all $I_k$ are Mumford-Tate generic, and 
the hyperk\"ahler structure $\c H_i$
induces the complex structures $I_i$, $I_{i-1}$.

\endproof


\section{Category of reflexive sheaves}
\label{_refle_Section_}


In this Section
we relate several classic results dealing
with reflexive sheaves and their singularities.
The proofs and further details can be found e.g. in 
\cite{_OSS_}.

\hfill

\definition\label{_refle_Definition_}
Let $X$ be a complex manifold, and $F$ a coherent sheaf on $X$.
Consider the sheaf $F^*:= \c Hom_{\calo_X}(F, \calo_X)$.
There is a natural functorial map 
$\rho_F:\; F \arrow F^{**}$. The sheaf $F^{**}$
is called {\bf a reflexive hull}, or {\bf 
reflexization},
of $F$. The sheaf $F$ is called {\bf reflexive} if the map
$\rho_F:\; F \arrow F^{**}$ is an isomorphism. 

\hfill

\remark
For all coherent sheaves $F$, the map
$\rho_{F^*}:\; F^* \arrow F^{***}$ is an isomorphism
(\cite{_OSS_}, Ch. II, the proof of Lemma 1.1.12).
Therefore, a 
reflexive hull of a sheaf is always 
reflexive.

\hfill

Reflexive hull can be obtained by 
restricting to an open subset and taking the
pushforward.

\hfill

\lemma\label{_refle_pushfor_Lemma_}
Let $X$ be a complex manifold, $F$ a coherent sheaf on $X$,
$Z$ a closed analytic subvariety, $\codim Z\geq 2$, and
$j:\; (X\backslash Z) \hookrightarrow X$ the natural
embedding. Assume that the pullback $j^* F$ is
reflexive on $(X\backslash Z)$. Then the pushforward
$j_* j^* F$ is also reflexive. 

\hfill

{\bf Proof:} This is \cite{_OSS_}, Ch. II, Lemma 1.1.12.
\endproof

\hfill

\definition
Let $F$ be a coherent sheaf on a complex manifold.
For any analytic subvariety $Z\subset X$, 
denote by $j:\; (X\backslash Z) \hookrightarrow X$ the natural
embedding. Consider the standard sheaf morphism
$F \stackrel \phi \arrow j_* j^* F$. The sheaf
$F$ is called {\bf normal} if $\phi$ is an isomorphism, for all
$\codim Z\geq 2$.

\hfill

From \ref{_refle_pushfor_Lemma_}, the following 
statement is clear.

\hfill

\lemma\label{_normal_shea_Lemma_}
(\cite{_OSS_})
Let $F$ be a coherent sheaf on a complex manifold.
Then $F$ is normal if and only if $F$ is reflexive.

\endproof

\hfill

This leads to the following definition

\hfill

\definition \label{_refle_category_Definition_}
Let $X$ be a complex manifold without subvarieties of codimension 1.
The category of reflexive sheaves on $X$ is defined as
\[ 
\c C^r(X) := \lim_\arrow \c C(X \backslash Z),
\]
where $\c C(X \backslash Z)$ is the category of coherent sheaves on
$X\backslash Z$ which can be extended to coherent 
sheaves on $X$,
and the limit is taken over all
increasing sequences of closed subvarieties
\[ Z \subset X, \ \ \ \codim _X Z \geq 2.
\]

\hfill

\lemma\label{_singu_refle_codim_3_Lemma_}
Let $F$ be a reflexive sheaf on $M$, and $X$ its singular set.
Then $\codim_M X\geq 3$

\hfill

{\bf Proof:} This is \cite{_OSS_}, Ch. II, 1.1.10. \endproof


\section{Twistor formalism and cohomology}
\label{_holo_and_twi_Section_}


\subsection{Twistor spaces}
\label{_twi_spa_Subsection_}

Let $M$ be a hyperk\"ahler manifold.  We identify the 2-dimensional 
sphere $S^2$ with the set of all quaternions $J$ with $J^2 = -1$. 
Consider the product manifold $X = M \times S^2$.
For every point
$x = m \times J \in X = M \times S^2$ the tangent space $T_xX$ is
canonically decomposed $T_xX = T_mM \oplus T_JS^2$. Identifying $S^2$
with $\C P^1$ in a standard fashion, 
we may assume that $S^2$ is equipped
with a natural complex structure.
Let $I_{S^2}:\; T_JS^2 \arrow T_JS^2$ be the complex structure operator.
Let $I_M:\; T_mM \arrow T_mM$ be the complex structure on $M$ induced by 
$J \in S^2 \subset {\Bbb H}$.

The operator $I_x = I_M \oplus I_J:\; T_MX \arrow 
T_xX$ satisfies $I_x \circ I_x = -1$. 
It depends smoothly on the point $x$, hence defines an almost complex
structure on $X$. This almost complex structure is known to be integrable
(see \cite{_Salamon_}). 

\hfill

\definition\label{_twistor_Definition_}
The complex manifold $(X, I_x)$ is called {\it the twistor
space} for the hyperk\"ahler manifold $M$, denoted by $\Tw(M)$.
This manifold is equipped with a real analytic projection
$\sigma:\; \Tw(M)\arrow M$ and a complex analytic
projection $\pi:\; \Tw(M) \arrow \C P^1$.

\hfill

For any $I\in \C P^1$, consider the complex submanifold
$\pi^{-1}(I)\subset \Tw(M)$. Clearly, $\pi^{-1}(I)$
is naturally isomorphic to $(M,I)$. Further on, 
we shall consider $(M,I)$ as a submanifold in $\Tw(M)$.

\subsection{Twistor formalism for vector bundles}
\label{_twi_fo_for_bu_Subsection_}

Let $M$ be a hyperk\"ahler manifold, and $\Tw(M)= M\times \C P^1$
its twistor space, equipped with a natural 
(non-holomorphic) projection $\sigma:\; \Tw(M) \arrow M$. 
Given a bundle $B$ with a connection $\nabla$, we can 
lift $\nabla$ to the pullback $\sigma^* B$ to obtain
a vector bundle $(\sigma^* B, \sigma^* \nabla)$
with connection. Take a $(0,1)$-part of $\sigma^* \nabla$.
By a Newlander-Nirenberg theorem (\ref{_Newle_Nie_for_bu_Theorem_}), 
the operator $(\sigma^* \nabla)^{0,1}$
defines a holomorphic structure on $\sigma^* B$
if and only if $((\sigma^* \nabla)^{0,1})^2 =0$.

In \cite{_NHYM_}, it is shown that this is equivalent to 
$SU(2)$-invariance of the curvature of the connection $\nabla$.

We obtain a functor (twistor transform)
\begin{equation}\label{_twi_functor_Equation_}
\sigma^*:\; (B, \nabla) \arrow (\sigma^* B, (\sigma^* \nabla)^{1,0}).
\end{equation}
It turns out that this functor is invertible: the connection
$\nabla$ can be reconstructed from the holomorphic structure
operator $(\sigma^* \nabla)^{0,1}$.
The following theorem was proven in \cite{_NHYM_}.

\hfill

\theorem\label{_twi_functor_Theorem_}
The functor \eqref{_twi_functor_Equation_}
gives an equivalence of the following categories
\begin{description}
\item[(i)] Vector bundles $B$ on $M$ with connection $\nabla$
whose curvature \[ \Theta\in \Lambda^2(M)\otimes \End B\]
is an $SU(2)$-invariant 2-form on $M$.

\item [(ii)] Holomorphic vector bundles $H$ on $\Tw(M)$
such that for all $x\in M$, the restriction of $H$ to a
holomorphic curve $\sigma^{-1}(x)\cong \C P^1 \subset \Tw(M)$
is trivial, as a holomorphic vector bundle:
\[ H\restrict{\sigma^{-1}(x)} \cong \calo_{\C P^1}^{\oplus n}.
\]
\end{description}

This equivalence is called {\bf the twistor transform}.

\endproof

\hfill

Assume now that $M$ is compact, and let $I$ be an induced complex structure.
Given a polystable bundle $B$ on $(M,I)$, with $c_1(B)$, $c_2(B)$
$SU(2)$-invariant, we can apply \ref{_UY_Theorem_} 
and obtain a unique Yang-Mills connection $\nabla$
on $B$. The curvature of $\nabla$ is 
$SU(2)$-invariant by 
\ref{_inva_then_hyperho_Theorem_}. Consider its twistor transform
\[
  \Tw(B):= (\sigma^* B, (\sigma^* \nabla)^{0,1})
\]
which is defined as in \ref{_twi_functor_Theorem_}.
We obtained a functor $B \arrow \Tw(B)$
from the category of hyperholomorphic bundles on $(M,I)$
(\ref{_hh_polysta_Definition_}) 
to the category of holomorphic
bundles on $\Tw(M)$.

\subsection{Cohomology of $\Tw(B)$}

Let $M$ be a compact hyperk\"ahler manifold, 
$I$ an induced complex structure, $B$ a hyperholomorphic bundle
on $(M,I)$, and $\Tw(B)$ the corresponding bundle on a twistor
space. Further on in this paper, we construct a holomorphic
bundle on a twistor space starting from a semistable bundle
on $(M,I)$. For this purpose, we need to relate the cohomology of $\Tw(B)$
and the cohomology of $B$.

\hfill

\proposition\label{_direct_ima_twi_Proposition_}
Let $M$ be a compact hyperk\"ahler manifold,
\[ \pi:\; \Tw(M) \arrow \C P^1\] its twistor space,
$L$ an induced complex structure, and $B$ a
hyperholomorphic bundle on $(M,L)$. Consider
the corresponding bundle $\Tw(B)$ on $\Tw(M)$.
Let $R^i \pi_* (\Tw(B))$ be the higher direct image
of $\Tw(B)$, considered as a coherent sheaf on
$\C P^1$. Then
\[ R^i \pi_* \Tw(B) \cong \calo _{\C P^1}(i) \otimes_\C H^i(B),
\]
where  $H^i(B)$ is the $i$-th cohomology
of $B$ on $(M,L)$.

\hfill

{\bf Proof:}
Essentially, \ref{_direct_ima_twi_Proposition_}
is proven in \cite{_Verbitsky:Hyperholo_bundles_}.

Consider an induced complex structure $L'$.
Restricting $\Tw(B)$ to \[ (M,L')\subset \Tw(M),\] we obtain
a holomorphic vector bundle $B_{L'}$. Denote
its cohomology by $H^i(B_{L'})$. In
\cite{_Verbitsky:Hyperholo_bundles_} it was proven
that $H^i(B_{L'})$ is (non-canonically) 
isomorphic to $H^i(B)$. This implies
that $R^i \pi_* \Tw(B)$ is a holomorphic vector
bundle, of rank $\dim H^i(B)$.
By Grothendieck's theorem, $R^i \pi_* \Tw(B)$
is a direct sum of several copies of $\calo(i)$, for
various $i$. We need only to show that $R^i \pi_* \Tw(B)$
is a direct sum of several copies of $\calo(i)$.

To prove this,
we need to repeat the argument of
\cite{_Verbitsky:Hyperholo_bundles_},
taking care of the dependency of $H^i(B_{L'})$
from the induced complex structure $L'$.

We say that an $SU(2)$-representation $V$ has weight 
$i$ if $V$ is generated by the highest weight
vectors of weight $i$.  An element $v\in V$ has weight $i$
if $v$ belongs to a subrepresentation of weight $i$.

Consider the space $\Lambda^i_+(M)\subset \Lambda^i(M)$
consisting of all $i$-forms which have weight $i$
(by multiplicativity of weight, 
an $i$-form has weight $\leq i$).
The space $\Lambda^i_+(M)$ is compatible
with the Hodge grading, that is,
it admits the Hodge decomposition
for every induced complex structure $I$
\begin{equation}\label{_Lambda_+_decompo_Equation_}
\Lambda^i_+(M) = \oplus_{p+q=i}\Lambda^{p,q}_{I,+}(M).
\end{equation}
In \cite{_V:projective_}, the space $\Lambda^i_+(M)\otimes B$
was studied in great detail. We have shown that
$\Lambda^i_+(M)\otimes B$ admits a natural $SU(2)$-invariant
Laplacian $\Delta_+$, whenever $B$ is
equipped with a hyperholomorphic
connection. Moreover, $\Delta_+$
is compatible with the Hodge decomposition,
and for all induced complex structures 
$I$ we have a canonical identification
\begin{equation}\label{_H^i_via_harmo_Equation_}
\c H^{0,p}_{+,I}(M,B) = H^p(B_I)
\end{equation}
where $\c H^{0,p}_{+,I}(M,B)$ is the 
space of $\Delta_+$-harmonic forms of type $(0,p)$,
and $H^p(B_I)$ the space of 
cohomology of $B$ considered as a holomorphic
bundle on $(M,I)$.

Let \[ \c W:= \c H^p_+(M,B)\otimes_\C \calo_{\C P^1}\]
be the trivial vector bundle with the fiber $\c H^p_+(M,B)$.
Using the Hodge decomposition 
\eqref{_Lambda_+_decompo_Equation_}
on the space of $\Delta_+$-harmonic forms
$\c H^p_+(M,B)$, 
we can decompose $\c H^p_+(M,B)\otimes_\C \calo_{\C P^1}$
onto its Hodge components\footnote{This decomposition
is not holomorphic. In fact, 
the Hodge decomposition
\eqref{_VHS_on_H^p(B)_Equation_} defines a
non-polarized variation of Hodge structures on 
$\c H^p_+(M,B)\otimes_\C \calo_{\C P^1}$.}
\begin{equation}\label{_VHS_on_H^p(B)_Equation_}
\c W\restrict I
= \oplus_{p+q=i}\c H^{p,q}_{+,I}(M,B) 
\end{equation}
By \eqref{_H^i_via_harmo_Equation_}, 
the bundle $R^p \pi_* \Tw(B)$
is identified with the $(0,p)$-part $\c W^{0,p}$ of the
trivial vector bundle $\c W= \c H^p_+(M,B)\otimes_\C \calo_{\C P^1}$.

Consider the natural embedding $\Xi$ of $\c H^{0,p}_{+,I}(M,B)$
to the space \[ \Lambda_\pi^{0,p}(\Tw (M))\otimes \Tw(B)\]
of relative $\Tw(B)$-valued $(0,p)$-forms on the space $\Tw(M)$.
On $(0,p)$-forms, the $SU(2)$-\-inva\-riant
Laplace operator $\Delta_+$
is equal to the usual Dolbeault Laplacian
\[ \Delta_{\bar \6}= \6 \6^* +\6^* \6.\] Therefore,
$\Xi$ is a holomorphic embedding, and we can reconstruct
the holomorphic structure on $R^p \pi_* \Tw(B)$
from that on $\Lambda_\pi^{0,p}(\Tw (M))\otimes \Tw(B)$.

Since the holomorphic structure on 
$R^p\pi_* \Tw(B) = \c W^{0,p}$
is compatible with the standard embedding
\[ \Xi:\; \c H^{0,p}_{+,I}(M,B) \arrow 
   \Lambda_\pi^{0,p}(\Tw (M))\otimes \Tw(B),
\]
the holomorphic structure on $R^p\pi_* \Tw(B)= \c W^{0,p}$
is inherited from $\c W$, as follows.
Take a section $Y$ of $\c W^{0,p}$,
and let $\bar\6_{\c W} Y$, $\bar\6_{\c W^{0,p}} Y$ be the 
holomorphic structure operator on $\c W$, $\c W^{0,p}$ applied to 
$Y$. Then $\bar\6_{\c W^{0,p}} Y$
is equal to the orthogonal projection of 
\begin{equation}\label{_holo_acti_in_W^0,p_Equation_}
 \bar\6_{\c W} Y
  \subset \Lambda^{0,1}\C P^1\otimes \c W
\end{equation}
to $\c W^{0,p}\subset \c W$.

Using the natural $SU(2)$-action on $\C P^1$
and $\c H^p_+(M,B)$, we may consider $\c W$
as an $SU(2)$-equivariant bundle. 
Let $I$ be an arbitrary induced complex structure,
and denote by $\c I$ the Lie algebra element 
associated with $I$\footnote{The action of
$\c I$ induces the Hodge decomposition
associated with $I$ in a standard way.}. Then 
$\c W^{0,p}\restrict I$ is by definition 
the space of all vectors $v\in \c W\restrict I$
which satisfy $\c I (v) = - p\1 v$.

We arrive at the following situation.
Let $W$ be an irreducible $SU(2)$-re\-pre\-sen\-tation of weight $p$, 
and
\[ 
   \c W = W \otimes_\C \calo_{\C P^1}
\] 
a trivial
bundle with an $SU(2)$-equivariant action,
induced by the $SU(2)$-action on $W$ and
$\C P^1$. Consider a $C^\infty$-subbundle
$\c W^{0,p}\subset \c W$ 
of all vectors $v\in \c W_\restrict I$
which satisfy 
\begin{equation}\label{_c_I_acts_on_W^0,p_Equation_}
\c I (v) = - p\1 v,
\end{equation}
and let the holomorphic structure
on $\c W^{0,p}$ act as in 
\eqref{_holo_acti_in_W^0,p_Equation_}.
Then \ref{_direct_ima_twi_Proposition_} 
is implied by the following
claim, which is a version of Borel-Bott-Weyl theorem
for the group $SL(2)$.

\hfill

\claim\label{_bu_calo(p)_lin_alge_Claim_}
In the above assumptions, the bundle
$\c W^{0,p}$ is isomorphic to a 
sum of several copies of $\calo(p)$.

\hfill

{\bf Proof:}
Clearly, it suffices to prove 
\ref{_bu_calo(p)_lin_alge_Claim_} when $W$ is 
an irreducible representation of weight $p$. In this case,
$\c W^{0,p}$ is obviously a line bundle.

Assume that $V$ is an irreducible representation of weight 1,
and let $\c V^{0,1}$ be the line bundle constructed
as above from $V$.
{}From \eqref{_c_I_acts_on_W^0,p_Equation_}
one can easily see that 
\begin{equation}\label{_p-th_power_Equation_}
   \c W^{0,p} \cong (V^{0,1})^{\otimes p}. 
\end{equation}
By \eqref{_p-th_power_Equation_},
it suffices to prove 
\ref{_bu_calo(p)_lin_alge_Claim_} for $p=1$,
and $W$ 2-dimensional.

Consider $\C P^1$ as a space of 1-dimensional
subspaces in $W$. Using the action of $SU(2)$
on $W$, we equip the bunle $\calo(1)$ 
of linear forms on these subspaces with a natural 
$SU(2)$-equivariant structure. Consider
a point $(1,0)\in \C P^1$ corresponding
to an induced complex structure $I$.
By definition, the Lie algebra element 
$\cal I$ acts on the line $(\C, 0)$ as $\1$,
hence it acts on its dual space
as $-\1$. Since a $G$-equivariant bundle 
over a manifold with a 
homogeneous $G$-action is uniquely 
determined by the action of 
the stabilizer of a point, 
we obtain that $\calo(1)$ is isomorphic to 
$\c W^{0,1}$ as an equivariant bundle. 
Therefore, $c_1(W^{0,1}) = c_1(\calo(1))=1$, 
and these line bundles are isomorphic
as holomorphic bundles as well.
\ref{_bu_calo(p)_lin_alge_Claim_} is proven.
We proved \ref{_direct_ima_twi_Proposition_}.
\endproof


\section{Local geometry of the twistor space}
\label{_twi_local_Section_}


\subsection{Twistor lines and local geometry}
\label{_twi_spa_local_Subsection_}

Consider the standard holomorphic projection
$\pi:\; \Tw(M) \arrow \C P^1$. 
To study $\Tw(M)$, we use the space of holomorphic 
sections \[ s:\; \C P^1 \arrow \Tw(M)\]
of $\pi$. Given $x\in M$, we have a tautological
section $s_x$ mapping $I\in \C P^1$ to 
$(x, I) \in \Tw(M) = M \times \C P^1$.
Such sections are called {\bf the horizontal
sections of $\pi$}, or {\bf horizontal twistor lines}. 
A normal bundle to $s$ is
isomorphic to $\calo(1)^{\oplus n}$
(\cite{_HKLR_}). By deformation theory,
this implies the following.

\begin{equation}\label{_twi_lines_geometry_minipage_Equation_}
\begin{minipage}[t]{0.7\linewidth}
Given two distinct points $I_1$, $I_2$ on $\C P^1$,
and two sufficiently close points $x_1$, $x_2\in M$,
there exists a unique section $s:\; \C P^1 \arrow \Tw(M)$
which is close to a horizontal section and
passes through $(I_1, x_1)$ and $(I_2, x_2)$
(see \cite{_Verbitsky:hypercomple_}
for a more rigorous version of this statement).
\end{minipage}
\end{equation}

\hfill

On an intuitive level, the twistor space is very close to a
projective space; indeed, for $M = \Bbb H$, $\Tw(M)$
is isomorphic to $\C P^3$ without a line.

\hfill

Fix a point $x_0\in M$, and the correspoinding
point in twistor space $(I, x_0)\in \Tw(M)$.
We shall study the twistor sections passing through
$(I, x_0)\in \Tw(M)$ in a neighbourhood of 
the horisontal section $s_{x_0}$.

For every $L, L'\in \C P^1$, let $\calo_{x_0, L}$,
$\calo_{x_0, L'}$ be the ring of germs of complex analytic
functions on $(M, L)$ and $(M, L')$ in $x_0$. In 
\cite{_Verbitsky:hypercomple_} we constructed a natural isomorphism
\begin{equation}\label{_Psi_on_germs_Equation_}
\Psi_I(L, L'):\; \calo_{x_0, L} \arrow \calo_{x_0, L'}
\end{equation}
associated with the local map 
$\tilde \Psi_I:\; (M,L') \arrow (M,L)$
defined in an appropriate neighbourhood of $x_0$
as follows.

Consider a point $(m, L')$ close to $(x_0, L')$.
By \eqref{_twi_lines_geometry_minipage_Equation_},
there exists a unique twistor section
$s:\; \C P^1 \arrow \Tw(M)$ 
passing through $(m, L')$ and $(x_0, I)$
and close to $s_{x_0}= \{x_0\} \times \C P^1$.
We define $\tilde \Psi_I(m):= s(L) \in (M, L)$
(see Figure 1). This map is holomorphic
and invertible in a sufficiently 
small neighbourhood of $(x_0, L')$
(q.v. \cite{_Verbitsky:hypercomple_}).

\begin{figure}[htb]
\centerline{\epsfig{file=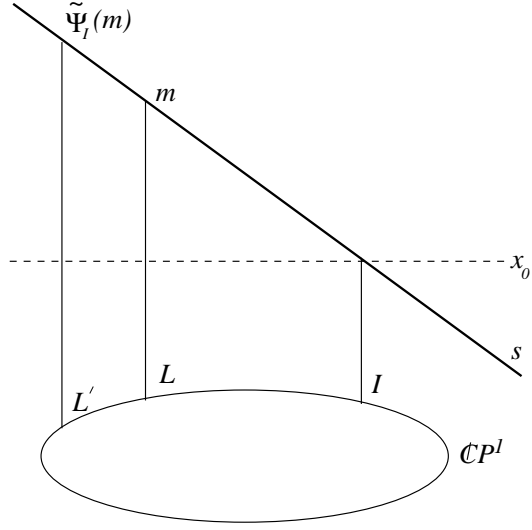,width=7cm}}
\caption{Constructing $\tilde \Psi_I$ with twistor sections}
\end{figure}

Let now $B$ be a bundle on $\Tw(M)$ obtained from the twistor
transform (\ref{_twi_functor_Theorem_}). This means that
the restriction of $B$ to any horizontal twistor
section $s_x$ is trivial. Clearly, a small deformation
of a trivial bundle on $\C P^1$ is again trivial.
Therefore, the restriction of $B$ to a twistor
section close to $s_x$ is trivial as well. 

This allows one to extend the map
$\Psi_I(L, L'):\; \calo_{x_0, L} \arrow \calo_{x_0, L'}$
to the space of germs of holomorphic sections of $B$
(see \cite{_V:Hyperholo_sheaves_} for 
a more detailed construction). We obtain an isomorphism
\begin{equation}\label{_Psi_on_germs_of_E_Equation_}
\Psi_I(L, L', B):\; B_{x_0, L}
  \arrow B_{x_0, L'}, 
\end{equation}
where $B_{x_0, L}$ and $B_{x_0, L'}$ are spaces
of germs of $B\restrict{(M,L)}$ and $B\restrict{(M,L')}$:
\[
B_{x_0, L} = \calo_{x_0, L}\otimes B\restrict{(M,L)} \ \ \ 
B_{x_0, L'} = \calo_{x_0, L'}\otimes B\restrict{(M,L')}.
\]
The isomorphisms $\Psi_I(L, L')$ and
$\Psi_I(L, L', B)$ depend on $L, L'\in \C P^1$ 
holomorphically (q. v. \cite{_V:Hyperholo_sheaves_}).

Consider the ring
$\calo_{\left(\Tw(M)\backslash (M, I)\right),s_{x_0}}$
of germs of $\calo_{\left(\Tw(M)\backslash (M, I)\right)}$
in a neighbourhood of 
\[ 
   s_{x_0}\backslash {(I, x_0)}\subset \Tw(M)\backslash (M, I).
\]
Denote the infinitesimal neighbourhood
of $s_{x_0}\backslash {(I, x_0)}$ in $\Tw(M)\backslash (M, I)$
by 
\[ 
   \Tw(M)_{x_0, I} = \Spec(\calo_{\Tw(M)\backslash (M, I),s_{x_0}}).
\]
This space is fibered over $\C P^1\backslash I$
with fibers isomorphic to a germ of a smooth
complex manifold. 

The maps \eqref{_Psi_on_germs_Equation_}
produce a canonical 
trivialization of $\Tw(M)_{x_0, I}$
over $\C P^1\backslash I$. Denote the trivialization
map by 
\begin{equation}\label{_Phi_triviali_Equation_}
  \Phi:\; \Tw(M)_{x_0, I} \arrow S\times (\C P^1\backslash I)
\end{equation}
where $S$ is a germ of a complex 
manifold isomorphic to any of the fibers of the
standard projection
\[ \pi:\; \Tw(M)_{x_0, I} \arrow \C P^1 \backslash I.
\]
Let $\xi:\;  \Tw(M)_{x_0, I} \arrow S$
be the composition of $\Phi$ and the projection
$S\times (\C P^1\backslash I)\arrow S$.
Using the maps $\Psi_I(L, L', B)$
of \eqref{_Psi_on_germs_of_E_Equation_}, we obtain
a trivialization of $B\restrict {\Tw(M)_{x_0, I}}$
over $\C P^1 \backslash I$. More precisely,
we obtain a bundle $B_S$ over $S$, and
a natural isomorphism
\begin{equation}\label{_E_loca_triv_Equation_}
B\restrict {\Tw(M)_{x_0, I}}\cong \xi^*B_S.
\end{equation}

\subsection{Local trivialization of the sheaf of differential operators}

Let $M$ be a hyperk\"ahler manifold, $\dim_\C M =2n$, 
$\pi:\; \Tw(M) \arrow \C P^1$ its twistor space,
$x\in M$ a point, $s_x:= \{x\}\times \C P^1\subset \Tw(M)$
the corresponding rational curve in the twistor space,
and $B$ a bundle on $\Tw(M)$ obtained 
from the twistor transform. By definition,
this means that the restriction
of $B$ to a curve $s_{x'}\subset \Tw(M)$ is a trivial 
bundle, for all $x'\in M$.\footnote{Any successive extension
of the bundles of form $\Tw(B_i)$, where $B_i$ are
hyperholomorphic, satisfies this condition.}
 Fix an induced complex
structure $I$. In Subsection \ref{_twi_spa_local_Subsection_},
we have constructed a trivialization of $B\restrict {\Tw(M)\backslash (M,I)}$
in a neighbourhood of $s_x\backslash I$.
In this Subsection we study the asymptotic
properties of this trivialization.

Let $\calo_{s_x}$ be the ring of germs of $\calo_{\Tw(M)}$
in a neighbourhood of $s_x$. Denote by $\goth{m}_{s_x}$
the ideal of all functions vanishing in $s_x$, and let
$\hat \calo_{s_x}$ be the $\goth{m}_{s_x}$-completion
of $\calo_{s_x}$.

We consider $\hat \calo_{s_x}$ with the 
natural $\goth{m}_{s_x}$-adic topology on it.
Clearly, we have 
\begin{equation}\label{_m/m^2_conormal_Equation_}
   \goth{m}_{s_x}/\goth{m}_{s_x}^2
   \cong \Omega^1_\pi \Tw(M)\restrict{s_x}.
\end{equation}
Using e.g. \ref{_bu_calo(p)_lin_alge_Claim_}, it is easy to check
that the sheaf of conormal $(1,0)$-vectors to $s_x$
is $\calo(-1)^{\oplus n}$:
\begin{equation}\label{_Omega_calo_-1_Equation_}
\Omega^1_\pi \Tw(M)\restrict{s_x} \cong \calo(-1)^{\oplus n}.
\end{equation}
Taking a symmetric power of \eqref{_Omega_calo_-1_Equation_},
and using \eqref{_m/m^2_conormal_Equation_}, we obtain
\begin{equation} \label{_maxi_ideal_on_CP^1_Equation_}
 \goth{m}_{s_x}^i/\goth{m}_{s_x}^{i+1} \cong \calo(-i)\otimes_\C S^i (\C^n),
\end{equation}
where $S^i (\C^n)$ denotes the $i$-th symmetric power of
the $n$-dimensional complex space $\C^n$.

Consider now the sheaf $\hat \calo_{s_x}^\star$
of continuous $\calo_{\C P^1}$-linear homomorphisms from $\hat \calo_{s_x}$
to $\calo_{\C P^1}$. The sheaf $\hat \calo_{s_x}^\star$
has a natural $\calo_{\Tw(M)}$-structure; for any section
$\gamma\in\hat \calo_{s_x}^\star$, $\delta \in \calo_{\Tw(M)}$,
$\epsilon \in \hat \calo_{s_x}$
we write 
\[ \delta \cdot \gamma (\epsilon) = \gamma(\delta\epsilon).
\]
This gives a structure of quasicoherent sheaf 
on $\hat \calo_{s_x}^\star$.

We consider
$\hat \calo_{s_x}^\star$ as an infinite-dimensional vector bundle
on $\C P^1$. The trivialization of $\calo_{s_x}$
outside $I$ gives a trivialization of
$\hat \calo_{s_x}$ outside $I$. However, 
the bundle $\hat \calo_{s_x}$ is
very negative by \eqref{_maxi_ideal_on_CP^1_Equation_},
and therefore admits no non-trivial sections.
By the same reason, the bundle $\hat \calo_{s_x}^\star$
is very positive, hence 
it admits many sections. One should expect that 
any section of $\hat \calo_{s_x}^\star\restrict{\C P^1\backslash I}$
compatible with a trivialization can be extended to
$\C P^1$. The aim of this section is to show that this
is indeed so.

\hfill

\definition
Given a bundle $F$ on $\C P^1$ with a trivialization
\[ 
  F\restrict{\C P^1\backslash I} 
  \cong \calo_{\C P^1\backslash I}\otimes_\C {\mathbf F}
\]
and a section $f\in F\restrict{\C P^1\backslash I}$,
we say that $f$ is {\bf compatible with the trivialization}
if there exists $f\in {\mathbf F}$ such that
\[ f= 1\otimes_\C f\in 
   \calo_{\C P^1\backslash I}\otimes_\C {\mathbf  F}.
\]

\hfill

Consider the sheaf 
\[ B_{s_x}^\star:= B^* \otimes_{\calo_{\Tw(M)}} \hat \calo_{s_x}^\star
\]
where $B$ is a bundle on $\Tw(M)$ considered in the beginning of this
Subsection, and $B^*$ its dual. Using the arguments of
Subsection \ref{_twi_spa_local_Subsection_}, we obtain a natural
trivialization of $B_{s_x}^\star$ over $\C P^1\backslash I$.
We study the asymptotical behaviour of this trivialization.

\hfill

\proposition\label{_B^star_trivialization_extended_Proposition_}
In the above assumptions, let $r\in  B_{s_x}^\star\restrict{\C P^1\backslash I}$
be a section which is compatible with the trivialization. Then $r$ can 
be extended to a section $\tilde r$ of $B_{s_x}^\star$.
Moreover, $\tilde r$ vanish at $I$.

\hfill

{\bf Proof:} To prove
\ref{_B^star_trivialization_extended_Proposition_},
we interpret $B_{s_x}^\star$ in terms of differential operators.

\hfill

\definition
Let $A$ be an algebra over 
the ring $B$, and $M$, $N$ $A$-modules. 
Following Grothendieck, we define 
the space $D_{A,B}(M, N)\subset \Hom_B(M,N)$ 
of $B$-linear differential operators from $M$ to $N$ 
inductively as follows.

An 0-th order differential operator $f\in D^0_{A,B} (M,N)$ is an
$A$-linear map from $M$ to $N$. An $n$-th order differential
operator $\delta:\; M \arrow N$ is a $B$-linear 
map such that for all $a\in A$,
the commutator $[a, \delta]$ is an $(n-1)$-st order differential
operator. The commutator $[a, \delta]:\; M \arrow N$
is defined as follows:
\[ [a, \delta] (m) = \delta (am) - a\delta(m).
\]

\hfill

The following observation is quite elementary.

\hfill

\claim \label{_diffe_ope_explici_Claim_}
In the assumptions of 
\ref{_B^star_trivialization_extended_Proposition_}, 
the sheaf $B_{s_x}^\star$ is naturally isomorphic
to the sheaf 
\[ D_{\calo_{\Tw(M), \calo_{\C P^1}}}(B, \calo_{s_x})
\]
of $\calo_{\C P^1}$-linear differential operators from
$B$ to the structure sheaf $\calo_{s_x}$ of $s_x\subset \Tw(M)$,
considered as a $\calo_{\Tw(M)}$-module. 

\hfill

{\bf Proof:} Denote by 
\[ D^n(B, \calo_{s_x})\subset D_{\calo_{\Tw(M), \calo_{\C P^1}}}(B,
  \calo_{s_x})
\]
the sheaf of $n$-th order differential operators. Since
\[ B_{s_x}^\star \cong \lim_\arrow 
    \Hom_{\calo_{\C P^1}}(B/\goth{m}_{s_x}^i B),
\]
to prove \ref{_diffe_ope_explici_Claim_}
it suffices to show that $D^n(B, \calo_{s_x})$
is dual over $\calo_{\C P^1}$ to $B/\goth{m}_{s_x}^{i+1} B$,
and that this duality is compatible with the corresponding 
arrows.

Take $\underline b\in B/\goth{m}^{n+1}_{s_x} B$,
and let $b\in B$ be its representative. Given
$\delta\in D^n(B, \calo_{s_x})$, we define
$\langle\underline b, \delta\rangle := \delta(b)$.
Since $\delta$ is of $n$-th order, $\langle\underline b, \delta\rangle$
is independent from the choice of $b\in B$.
Clearly, for all $\delta\neq 0$ there is $b$ such that
$\langle b,\delta\rangle \neq 0$. Comparing the dimensions of 
$D^n(B, \calo_{s_x})$ and $B/\goth{m}^{n+1}_{s_x} B$,
considered as a finite-dimensional bundles on $\C P^1$,
we find that this pairing is non-degenerate. We have shown
that $D^n(B, \calo_{s_x})$ is dual to $B/\goth{m}_{s_x}^{i+1} B$.
\ref{_diffe_ope_explici_Claim_} is proven. \endproof

\hfill

Now let $D\in D^n(B, \calo_{s_x})\restrict{\C P^1\backslash I}$
be a differential operator which is compatible with the
trivialization. We shall describe $D$ in explicit terms
as follows.

Let $L\neq I$ be an arbitrary induced complex structure,
and $U_L$ the germ of $(M,L)$ at $x\in M$.
For any $y\in U_L$, let $r_y:\; \C P^1 \arrow \Tw(M)$
be the unique line close to $s_x$ and passing
through $(y,L)$ and $x, I\in \Tw(M)$.
This gives a map $U_L\times \C P^1 \stackrel \Xi \arrow \Tw(M)$,
$(y, L') \arrow r_y(L')$ (see Figure 2).

\begin{figure}[htb]
\centerline{\epsfig{file=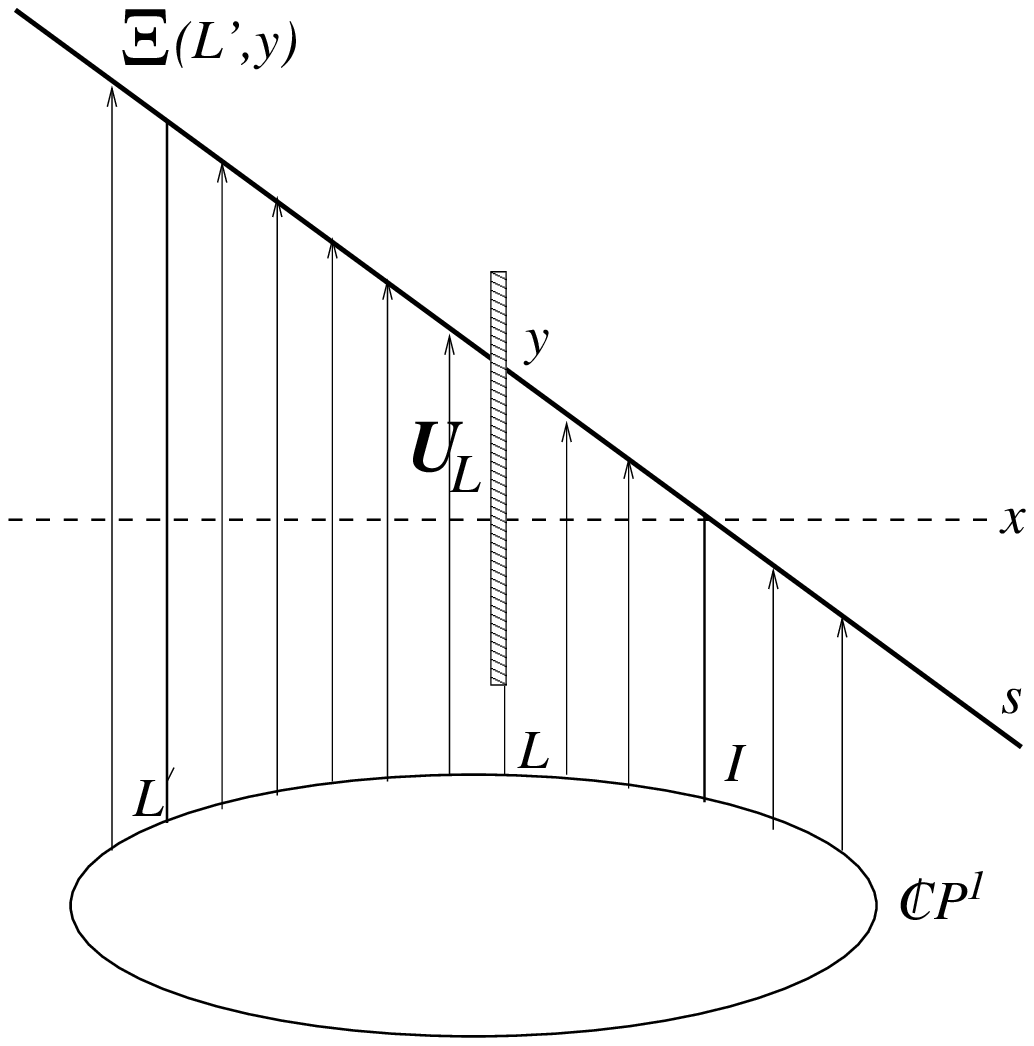,width=7cm}}
\caption{The map $U_L\times \C P^1 \stackrel \Xi \arrow \Tw(M)$}
\end{figure}

Given a differential operator 
$\underline{D_1}:\; B\restrict{U_L} \arrow \calo_{U_L}$,
consider the operator 
\[ D_1:= 
   \underline D_1 \boxtimes id :\; \Xi^*(B\restrict{U_L}) \arrow 
   \calo_{U_L\times \C P^1}
\]
We define the operator 
\[ \Xi(\underline D_1) \in D_{\calo_{\Tw(M), \calo_{\C P^1}}}(B, \calo_{s_x})
\]
as follows: given a section $b\in B$, we take the pullback
$\Xi^*(b)$, apply $D_1$ and then restrict to 
\[ s_x = \{x\}\times \C P^1\subset U_L \times \C P^1 = \Tw(U_L).\]
Tracing back to the definition, we obtain that the differential
operator 
\[ D\in D_{\calo_{\Tw(M), \calo_{\C P^1}}}(B, \calo_{s_x})\restrict{\C
    P^1\backslash I}, 
\]
is compatible with the trivialization
if and only if $D$ is constructed as above:
$D = \Xi(\underline D_1)\restrict{\C P^1\backslash I}$, where
where $\underline D_1$ is a differential
operator on $U_L$. When $D$ is a first order  differential
operator, that is, a vector field, this can be illustrated by the
following picture (Figure 3).

\begin{figure}[htb]
\centerline{\epsfig{file=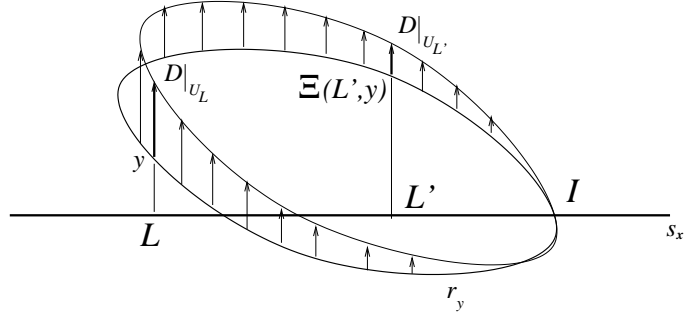,width=9cm}}
\caption{First order differential operator 
trivialized by local twistor geometry}
\end{figure}

In Figure 3, the little arrows represent the vertical
vector field corresponding to $D= \Xi(\underline D_1)$.
To differentiate a function $f$ along $D$, we look how $r_y(L')$
changes as $y$ goes along the vector field
$\underline D_1$ on $U_L$. Now 
\ref{_B^star_trivialization_extended_Proposition_}
is apparent: given a differential operator
$\underline D_1$ on $U_L$, and a germ of $B$
around $s_x\subset \Tw(M)$, we define  $\Xi(\underline D_1)(b)$ 
as above, and the differential operator $\Xi(\underline D_1)$
extends from $\C P^1\backslash I$ to $\C P^1$ 
in a natural way (in Figure 3, the vector fields
$D$ vanishes at
$I$). \ref{_B^star_trivialization_extended_Proposition_}
is proven. \endproof

\subsection{Twistor space and cohomology with support}

Given a coherent sheaf $F$ on $X$
and a finite set $A\subset X$, let
$\c H^i_A(F)$ be the sheaf of cohomology of $F$ with support in $A$.
Further on, we shall need the following characterization
of the cohomology with support. 

\hfill

\proposition\label{_Gro_support_Proposition_}
(Grothendieck's local duality)
Let $X$ be a compact complex manifold, $A\subset X$ 
a finite subset, and $B$ a holomorphic bundle on $X$.
Then $\c H^i_A(B)=0$ for all $i\neq \dim X$. 
For $i= \dim X$ the sheaf $\c H^i_A(B)$
can be expressed as follows. Take the
completion $\hat \calo_{X,A}$ of $\calo_X$ in $A$-adic topology.
Consider the dual space $\Hom_\C (\hat \calo_{X,A}, \C)$
of continuous $\C$-linear maps from $\hat \calo_{X,A}$ 
to $\C$, equipped with a natural $\calo_X$-action.
Let $D_A:=\Hom_\C (\hat \calo_{X,A}, \C)\otimes _{\calo_X}K_X$
be the product of $\Hom_\C (\hat \calo_{X,A}, \C)$
with the canonical sheaf. Then $\c H^i_A(B)$
is naturally isomorphic to $D_A$.

\hfill

{\bf Proof:} See \cite{_LC_}, \S 4, Example 3. \endproof

\hfill

From \ref{_Gro_support_Proposition_} we can easily obtain the following
corollary.

\hfill

\corollary\label{_open_coho_dim_less_n_Corollary_}
Let $X$ be a compact complex manifold, $F$ a coherent sheaf
on $X$ and $Z\subset X$ a non-empty complex analytic subset.
Then $H^n(X\backslash Z, F\restrict{X\backslash Z})=0$,
for $n\geq \dim_\C X$.

\hfill

{\bf Proof:} See \cite{_LC_}, Theorem 6.9. \endproof

\hfill

\remark
The statement of \ref{_open_coho_dim_less_n_Corollary_}
is often formulated as follows: one says that $X\backslash Z$
{\bf has cohomological dimension at most $\dim X-1$.}

\hfill

Let $M$ be a compact hyperk\"ahler manifold, $I$ an induced
complex structure, $A$ a finite set, and $B$ a hyperholomorphic
bundle on $M$. We study the cohomology of the holomorphic
bundle $\Tw(B)\restrict{\Tw(M\backslash A)}$. Consider $B$
as a holomorphic bundle on $(M,I)$.
Let $j:\; (M\backslash A)\hookrightarrow M$ be the natural
embedding. For an acyclic sheaf $\c F$, consider the exact sequence
\[ 0 \arrow H^0_A(\c F) \arrow F \arrow j^* \c F \arrow 0
\]
The corresponding long exact sequence gives
\begin{equation} \label{_j_H_a_long_exact_Equation_}
H^1_A(B) \arrow H^1(B) \arrow H^1(j^* B) \arrow H^2_A(B)
\end{equation}
For $\dim_\C M>2$, we have 
$H^2_A(B)= H^1_A(B)=0$ (\ref{_Gro_support_Proposition_}).
Then, \eqref{_j_H_a_long_exact_Equation_} 
implies $H^1(j^* B)\cong H^1(B)$.
For $\dim_\C M=2$, the situation is different.
We have a long  exact sequence
\[ 0\arrow H^1(B) \arrow H^1(j^* B) \arrow H^2_A(B) \arrow H^2(B) \arrow
   H^2(j^* B).
\]
However, the cohomological dimension of $M\backslash A$
is $\leq 1$ by \ref{_open_coho_dim_less_n_Corollary_},
hence $H^2(j^* B)=0$. We obtain the following result.

\hfill

\corollary
Let $X$ be a complex surface, $B$ a holomorphic
vector bundle, and $A$ a non-empty finite set.
Then we have an exact sequence
\begin{equation} \label{_j_H_a_long_exact_dim_2_Equation_}
0 \arrow H^1(B) \arrow H^1(j^* B) \arrow H^2_A(B) \arrow H^2(B)
\arrow 0.
\end{equation}
\endproof

\hfill

The main result of this Subsection is the following

\hfill

\proposition\label{_dire_on_tw_with_suppo_compu_Proposition_}
Let $M$ be a hyperk\"ahler manifold of real dimension 4, 
$\Tw(M) \stackrel \pi \arrow \C P^1$ its twistor space,
$I$ an induced
complex structure, $A\subset M$ a finite set, and $B$ a hyperholomorphic
bundle on $M$. Consider $B$ as a holomorphic bundle on $(M,I)$.
Let $\Tw(B)$ be the twistor transform of $B$. 
Denote by $\Tw_A(B)$ the restriction of $\Tw(B)$ 
to $\Tw(M\backslash A)$. Then 
\begin{equation}\label{_cohomo_vanishes_Equation_}
H^i(\Tw_A(B))=0 \text{\ \ \ for\ \ \ } i>1.
\end{equation}
Moreover, we have an
exact sequence of sheaves on $\C P^1$
\begin{equation}\label{_cohomo_w_suppo_on_twi_long_exa_Equation}
\begin{aligned}
0\arrow & H^1(B) \otimes_\C \calo(1) \arrow R^1 \pi_* \Tw_A(B) \\
{}& \arrow \bigoplus_{x_i\in A} \Hom_\C(\hat\calo_{s_{x_i}},\C)\otimes_{\calo_{\Tw(M)}}
  \Tw(B) \otimes_{\calo(\C P^1)} \calo(2) \\ 
{}& \arrow H^2(B) \otimes_\C
  \calo(2) \arrow 0  
\end{aligned}
\end{equation}
where $\hat\calo_{s_{x_i}}$ is the $s_{x_i}$-adic completion
of $\calo_{\Tw(M)}$, and $s_{x_i}\subset \Tw(M)$ is the rational
curve $x_i\times \C P^1 \subset \Tw(M)$, and
$\Hom_\C(\hat\calo_{s_{x_i}},\C)$ the sheaf of 
continuous (in adic topology) $\calo_{\C P^1}$-linear
homomorpisms from $\hat\calo_{s_{x_i}}$ to $\calo_{\C P^1}$.

\hfill

{\bf Proof:} The exact sequence 
\eqref{_cohomo_w_suppo_on_twi_long_exa_Equation}
is a relative version of \eqref{_j_H_a_long_exact_dim_2_Equation_}.
By \ref{_direct_ima_twi_Proposition_},
we have $R^i \pi_* \Tw(B) \cong H^i(B) \otimes_\C
  \calo(i)$; this gives us the rightmost and the
leftmost term of the exact sequence 
\eqref{_cohomo_w_suppo_on_twi_long_exa_Equation}.
To prove \eqref{_cohomo_w_suppo_on_twi_long_exa_Equation}
it remains to show that
\begin{equation} \label{_local_over_twisto_compu_Equation_}
\begin{aligned}
  \bigoplus_{x_i\in A} & \Hom_\C(\hat\calo_{s_{x_i}},\C)\otimes_{\calo_{\Tw(M)}}
  \Tw(B) \otimes_{\calo(\C P^1)} \calo(2)\\ \cong & \c H^1_{A\times \C
  P^1}(\Tw(B))
\end{aligned}
\end{equation}
The fiberwise canonical class of $\Tw(M)$ is $\calo(2)$
because the normal bundle to $s_{x_i}$ is $\calo(1)\oplus\calo(1)$.
Therefore, \eqref{_local_over_twisto_compu_Equation_}
is implied by \ref{_Gro_support_Proposition_}.

To prove that $H^2(\Tw_A B)=0$, notice that $R^2\pi_* \Tw_A(B)=0$
because the fibers of $\pi:\; \Tw(M\backslash A) \arrow \C P^1$
have cohomological dimension $\leq 1$
(\ref{_open_coho_dim_less_n_Corollary_}).
Therefore, 
\begin{equation}\label{_seco_coho_on_twi_via_suppo_Equation_} 
  H^2(\Tw_A(B)) = H^1(R^1\pi_* \Tw_A(B)).
\end{equation}
On the other hand, the sheaf $R^1\pi_* \Tw_A(B)$
can be obtained in an explicit way from the exact sequence 
\eqref{_cohomo_w_suppo_on_twi_long_exa_Equation}
as follows. 

Consider the $s_{x_i}$-adic filtration on $\hat\calo_{s_{x_i}}$.
The associated graded sheaf is isomorphic to
$\oplus_j {\goth m}_{s_i}^{j-1}/{\goth m}_{s_i}^j$,
where ${\goth m}_{s_i}$ is the maximal ideal of $s_{x_i}$.
However,
\[ 
  {\goth m}_{s_i}^{j-1}/{\goth m}_{s_i}^j \cong \calo(-j)\otimes_\C
  S^j \Lambda^{1,0}_{x_i}(M,I),
\]
because ${\goth m}_{s_i}^{j-1}/{\goth m}_{s_i}^j$
is the $j$-th symmetric power of 
\[ {\goth m}_{s_i}/{\goth m}_{s_i}^2\cong 
   \Lambda^{1,0}_\pi(\Tw(M))\restrict{s_{x_i}} \cong N^*(s_{x_i})
\]
and the normal bundle $N(s_{x_i})$ 
is isomorphic to $\calo(1)^{\oplus^2}$.

The ${\goth m}_{s_i}$-adic filtration on 
$\hat\calo_{s_{x_i}}$ induces a filtration on
\[ \Hom_\C(\hat\calo_{s_{x_i}},\C),\] with the associated graded
quotient sheaves isomorphic to 
\[ \calo(j)\otimes_\C  S^j T^{1,0}_{x_i}(M,I) \]
These bundles are positive, for all $j$. This allows us to represent
\[ \Hom_\C(\hat\calo_{s_{x_i}},\C)\] as a direct limit of
of positive bundles:
\[\Hom_\C(\hat\calo_{s_{x_i}},\C) = \lim_\arrow 
 \Hom_\C(\calo_{s_{x_i}}/{\goth m}_{s_i}^j,\C),
\]
with all $\Hom_\C(\calo_{s_{x_i}}/{\goth m}_{s_i}^j,\C)$
admitting a filtration with positive associated graded
quotients, hence positive.
We obtain that $\Hom_\C(\hat\calo_{s_{x_i}},\C)$
is a direct limit of positive bundles.
From the exact sequence
\eqref{_cohomo_w_suppo_on_twi_long_exa_Equation}
we obtain
\begin{equation} \label{_cohomo_w_suppo_positive_Equation_}
\begin{aligned}
0\arrow &H^1(B) \otimes_\C \calo(1) \arrow R^1 \pi_* \Tw_A(B) \\
\arrow & \oplus_{\alpha_k}\calo(\alpha_k)\arrow H^2(B) \otimes_\C
  \calo(2) \arrow 0,  \\ {} &\alpha_k \geq 2.
\end{aligned}
\end{equation}
This implies immediately that $R^1 \pi_* \Tw_A(B) $
is a direct sum of $\calo(k)$ with $k\geq 1$, and
its first cohomology vanish. Now 
\eqref{_seco_coho_on_twi_via_suppo_Equation_} 
implies $H^2(\Tw_A(B))=0$. We proved 
\ref{_dire_on_tw_with_suppo_compu_Proposition_}.
\endproof

\hfill

Let $M$ be a hyperk\"ahler manifold, and $x\in M$
the point. By $s_x\subset \Tw(M)$  we denote
the rational curve $\C P^1 \times \{ x\}$.
Consider the local cohomology sheaf 
\[ D_{s_{x}}:= \c H^{\dim_\C M}_{s_{x}}(\calo_{\Tw(M)})= 
   \Hom_\C(\hat\calo_{s_{x}},\C)\otimes_{\calo_{\Tw(M)}} K_\pi(\Tw(M))
\]
(the last equality holds by \ref{_Gro_support_Proposition_}). 
Using \eqref{_Phi_triviali_Equation_}, we can 
trivialise the infinite-dimensional bundle $\hat \calo_{s_{x}}$
outside of $I$ in a natural way. This gives
a trivialization of $D_{s_{x}}$ outside of $I$
(the bundle $K_\pi(\Tw(M))\cong \calo(\dim_\C M)$
is trivialized by taking the sections which have a zero
of order $n$ in $I$). Given a hyperholomorphic bundle
$B$ on $M$, we obtain a trivialization of $\Tw(B)$
in a neighbourhood of 
$s_x\backslash I \subset \Tw(M) \backslash (M,I)$
(see \eqref{_E_loca_triv_Equation_}). This gives a trivialization of
the term
\[ \bigoplus_{x_i\in A} \Hom_\C(\hat\calo_{s_{x_i}},\C)\otimes_{\calo_{\Tw(M)}}
  \Tw(B) \otimes_{\calo(\C P^1)} \calo(2)
\]
from the exact sequence 
\eqref{_cohomo_w_suppo_on_twi_long_exa_Equation}.
The sheaf $R^2 \pi_* \Tw(B) = H^2(B) \otimes_\C
  \calo(2)$ is trivialized outside $I$ by taking the sections
with second order zeroes at $I$. We obtained a trivialization
of the sheaf
\begin{equation}\label{_loca_coho_sheaf_trivialized_Equation_}
\begin{aligned}
\ker\bigg( 
\bigoplus_{x_i\in A} & \Hom_\C(\hat\calo_{s_{x_i}},\C)\otimes_{\calo_{\Tw(M)}}
  \Tw(B) \otimes_{\calo(\C P^1)} \calo(2) 
  \\ \arrow & R^2\pi_* \Tw(B) \bigg)\restrict{\C P^1\backslash I}.
\end{aligned}
\end{equation}

\hfill

\corollary
In the above assumptions, let $\nu$ be a section
of the sheaf \eqref{_loca_coho_sheaf_trivialized_Equation_}
which is compatible with the trivialization constructed above.
Then $\nu$ van be extended to a holomorphic section of the sheaf
\begin{equation*}
\begin{aligned}
\ker\bigg( 
\bigoplus_{x_i\in A} & \Hom_\C(\hat\calo_{s_{x_i}},\C)\otimes_{\calo_{\Tw(M)}}
  \Tw(B) \otimes_{\calo(\C P^1)} \calo(2) 
  \\ \arrow & R^2\pi_* \Tw(B) \bigg).
\end{aligned}
\end{equation*}
on $\C P^1$.

\hfill

{\bf Proof:} Follows immediately from
\ref{_B^star_trivialization_extended_Proposition_}. 
\endproof


\section{Coherent sheaves on generic compact tori}
\label{_tori_Section_}


\subsection{Reflexive sheaves on generic complex tori are bundles}
\label{_refle_tori_Subsection_}

\theorem\label{_cohe_on_toru_Theorem_}
Let $T$ be a compact  complex torus, $\dim_\C T >2$.
Assume that $H^{p,p}(T) \cap H^{2p}(T, \Z)=0$,
for all $0<p<\dim_\C T$.\footnote{This is equivalent 
to $T$ being Mumford-Tate generic.} 
Then all coherent sheaves on $T$ have isolated 
singularities and smooth reflexive hulls.

\hfill

{\bf Proof:} Consider a coherent sheaf $F$ on $T$. Since $T$ has
no non-trivial integer $(p,p)$-cycles, all subvarieties of $T$ are 
finite sets of points. Therefore, $F$ has isolated singularities. 
It remains to prove that the reflexive hull 
$E=F^{**}$ is smooth.
The following argument is well known
(see e.g. \cite{_Voisin_}).

\hfill

\lemma \label{_refle_she_on_torus_bundles_Lemma_}
Let $T$ be a compact complex 
torus with no non-trivial integer $(p,p)$-cycles,
for $0<p <3$,
and $E$ a reflexive coherent sheaf on $T$.
Then $E$ is a bundle.

\hfill

{\bf Proof:} We may assume that 
$\dim_\C \geq 3$.
Since $T$ has no non-trivial
integer $(p,p)$-cycles,
the first Chern class of all sheaves 
on $T$ vanishes. Therefore, all
coherent sheaves on $T$ 
are semistable.  Consider a filtration
\[
0 = E_0 \subset E_1 \subset ... \subset E_n = E
\]
with all the subfactors $E_i/E_{i-1}$ stable. 
First of all, we show that the reflexizations
of the sheaves $E_i/E_{i-1}$ are smooth, for all $i$. 
By construction, these sheaves are stable. 
Assume that $E$ is stable and reflexive.

In \cite{_Bando_Siu_}, Bando and Siu construct
a canonical Hermitian Yang-Mills  connection $\nabla$
on the non-singular part of any stable
reflexive sheaf on a compact K\"ahler
manifold $X$, with $L^2$-integrable
curvature $\Theta$. Such a connection
(called admissible Yang-Mills) is unique. 
In \cite{_Bando_Siu_} (see also \cite{_Tian:gauge_})
it was proven that the $L^1$-integrable form $\Tr(\Theta\wedge\Theta)$
(considered as a current on the K\"ahler manifold
$X$) is closed and represents the 
cohomology class 
\begin{equation}\label{_curva_Chern_Equation_} 
  [\Theta] = c_2(E) - \frac{r-1}{r} c_1^2(E)
\end{equation}
Consider the standard Hodge operator
$\Lambda$ on differential forms
(see \ref{_Yang-Mills_Definition_}).
Since $\nabla$ is Yang-Mills, we have
\begin{equation}\label{_Lubcke_Equation_} 
\Lambda^2 (\Theta^2) \geq 0,
\end{equation}
and the equality is reached only if $\Theta=0$.\footnote{This 
is the celebrated L\"ubcke inequality \cite{_Lubcke_},
which implies flatness
of stable bundles with zero Chern
classes and the Bogomolov-Miyaoka-Yau 
inequality; see  \cite{_Bando_Siu_}, \cite{_Tian:gauge_}}
Comparing \eqref{_curva_Chern_Equation_} 
and \eqref{_Lubcke_Equation_},
and using $c_1(E), c_2(E)=0$, we obtain
that $\Theta=0$, that is, $\nabla$ is flat.
Since $E$ is non-singular in codimension 2,
$\nabla$ has no local monodromy. Therefore,
$(E,\nabla)$ can be extended to a flat holomorphic
bundle on $T$. 
By \ref{_singu_refle_codim_3_Lemma_}, 
this extension is a reflexization of $E$.

We proved that $E^{**}$ is smooth, for all stable 
coherent sheaves $E$ on $T$. This implies that any
vector bundle $E$ on $T$ is filtered by subsheaves
$E_i$, with the quotient sheaves $E_i/E_{i-1}$ 
all having smooth reflexizations. Replacing $E_i$
by its reflexization $E_i^{**}\subset E$, we 
may assume that al $E_i$ are reflexive.
Using induction, we may assume also that $E_{n-1}$
is smooth. 

We have an exact sequence
\[
0 \arrow E_{n-1} \arrow E \arrow E_n / E_{n-1} \arrow 0
\]
with $E_{n-1}$ smooth, $E$ reflexive, and $F = E_n / E_{n-1}$
a stable sheaf having (as we have shown above) a smooth
reflexization.

Then $E$ is given by a class $\nu \in \Ext^1(F, E_{n-1})$.
Consider the exact sequence
\[ 0 \arrow F \arrow F^{**} \arrow C \arrow 0
\]
where $C$ is a torsion sheaf (cokernel of the reflexization map).
This gives a long exact sequence
\begin{equation}\label{_long_refle_on_torus_Equation_}
\Ext^1(F^{**}, E_{n-1}) \arrow \Ext^1(F, E_{n-1})
\stackrel \delta \arrow \Ext^2(C, E_{n-1}) 
\end{equation}
The kernel of $\delta$ in \eqref{_long_refle_on_torus_Equation_}
corresponds to all extensions 
$\gamma\in \Ext^1(F, E_{n-1})$
with a reflexization
isomorphic to an extension of $F^{**}$ with $E_{n-1}$.
Clearly, such extensions are reflexive only of $C=0$.
To prove that $C=0$, it suffices to show that $\delta(\nu)=0$.
However, $C$ is a torsion sheaf, therefore its support $S$
is a finite set. By Grothendieck's duality
(\ref{_Gro_support_Proposition_}), 
the group $\Ext^2(C, E_{n-1})\subset H^2_S(E_{n-1})$
vanishes, for $\dim T >2$. Therefore,
$E$ is smooth. This proves 
\ref{_refle_she_on_torus_bundles_Lemma_}.
\endproof

\hfill

The same argument also proves the following corollary.

\hfill

\corollary\label{_exa_on_to_bun_Corollary_}
In assumptions of \ref{_refle_she_on_torus_bundles_Lemma_},
let
\[
0 \arrow F_1 \arrow F_2 \arrow F_3 \arrow 0
\]
be an exact sequence of coherent sheaves, with $F_1$,
$F_2$ reflexive (hence, by \ref{_refle_she_on_torus_bundles_Lemma_},
locally free). Then the sheaf $F_3$ is also locally free.

\endproof

\hfill

Let $T$ be a compact complex 
torus of complex dimension $>2$ 
which has no non-trivial integer $(p,p)$-cycles, $0 <p <3$.
Since $c_1(B)=0$ for any sheaf $B$ on $T$, 
any sheaf on $T$ is semistable.  Consider the Jordan-Holder
filtration, with polystable associate quotient sheaves $B_i$.
By \ref{_refle_she_on_torus_bundles_Lemma_},
the reflexization $B_i^{**}$ is a bundle. However,
the sheaves $B_i$ are not necessarily bundles. 
Indeed, let $\Omega\subset \C^n$ be a Stein domain,
and consider an extension
\begin{equation}\label{_exte_non-tri_Equation}
   0 \arrow \calo_\Omega \arrow B \arrow \calo_\Omega(x) \arrow 0
\end{equation}
where $\calo_\Omega(x)$ is a kernel of a non-zero map
$\calo_\Omega \arrow k_x$ from $\calo_\Omega$ to a sky\-scraper sheaf.
If \eqref{_exte_non-tri_Equation} is non-trivial, 
then $B$ is a bundle, as one can easily check. 
It is easy to construct such an extension for
$n =2$. However, over a generic complex torus 
of dimension $>2$
this situation is impossible, as the following trivial 
proposition implies.

\hfill

\proposition\label{_filtra_w_line_bu_Proposition_}
Let $T$ be a $d$-dimensional compact complex 
torus, $d\geq 3$. Assume that $T$ 
has no non-trivial integer $(p,p)$-cycles:
\[ H^{p,p}(T) \cap H^{2p}(T, \Z) =0, \text{\ \ for\ \ } 0<p<d.
\]
Then any holomorphic vector bundle 
on $T$ admits a locally free filtration
by holomorphic bundles with associated 
graded sheaves locally free of rank 1.

\hfill

{\bf Proof:} Consider a Jordan-Holder filtration
\[ 0 = F_0\subset F_1 \subset ...\subset F_n= B,\]
with associated graded sheaves stable. If $F_i$ are not reflexive,
we replace $F_i$ by $F_i^{**}\subset B$.
Therefore, we may assume that the sheaves
$F_i$ are reflexive. \ref{_exa_on_to_bun_Corollary_}
implies that for all $i$, the quotient $F_i/F_{i-1}$
is a bundle. These bundles are also stable, and therefore
admit a Yang-Mills metrics. By L\"ubcke's 
(also Bogomolov's and Simpson's) argument (\cite{_Simpson_})
a Yang-Mills connection on a bundle with the
zero Chern classes $c_1$, $c_2$ is flat. Thus, it corresponds 
to an irreducible  unitary representation of $\pi_1(T)$. However, 
$\pi_1(T)$ is abelian, hence all its irreducible  
unitary representations are 1-dimensional.
We proved that $F_i/F_{i-1}$ is a line bundle.
\ref{_filtra_w_line_bu_Proposition_} is proven.
\endproof

\hfill

Consider a holomorphic vector 
bundle $B$ on a  compact  complex $d$-di\-men\-sional
torus $T$ with
\[ H^{p,p}(T) \cap H^{2p}(T, \Z) =0, \text{\ \ for\ \ } 0<p<d.
\]
In this Subsection we prove that $B$ admits a natural flat
connection $\nabla$ compatible with a holomorphic structure. 
If $\nabla$ is flat and Hermitian, then
it is Yang-Mills. In this case, 
$B$ is polystable by Uhlenbeck-Yau
theorem (\ref{_UY_Theorem_}).

Generally speaking, $B$ is not stable, and 
$\nabla$ is not Hermitian.

\hfill

Let $B_1$, $B_2$ be flat holomorphic Hermitian vector bundles
on $T$. Using the Hodge theory, we identify
$\Ext^1(B_1, B_2)$ with the space of harmonic 
$(0,1)$-forms with coefficients in $\Hom(B_1, B_2)$.
The maximum principle implies that any
harmonic $\Hom(B_1, B_2)$-valued form is parallel
with respect to the natural flat connection on
$\Lambda^1(M)\otimes B_1^* \otimes B_2$.
Let $B$ be the Yoneda extension of $B_2$ with $B_1$:
\begin{equation}\label{_exte_on_T_Equation_}
0\arrow B_1\arrow B\arrow B_2\arrow 0.
\end{equation}
The holomorphic structure operator in $B$ can be 
written explicitly as follows. Let $\nu\in \c H^{0,1}(\Hom(B_2,B_1)$
be the harmonic representative of the extension class
defining \eqref{_exte_on_T_Equation_}, and
$\bar\6_1$, $\bar\6_2$ the holomorphic structure
operators on $B_1$, $B_2$, 
\[ \bar\6_i :\; B_i \arrow B_i\otimes \Lambda^{0,1}(M).
\]
Consider the holomorphic structure operator
$\bar\6_{gr}:=\bar\6_1 + \bar\6_2$ on $B_{gr} := B_1\oplus B_2$.
We define an operator
\[ \bar\6:\; B_{gr} \arrow B_{gr} \otimes
 \Lambda^{0,1}(M)
\]
as
\[ \bar\6 := \bar\6_{gr} + \nu
\]
where 
\[ \nu\in \c H^{0,1}(\Hom(B_2,B_1))\subset \Lambda^{0,1}(\End(B_{gr}))\]
is understood as a connection form.
We have
\begin{equation}\label{_MC_Equation_} 
  \bar\6^2 = \bar\6_{gr}(\nu) + 2 \nu\wedge\nu =0
\end{equation}
(this is the famous Maurer-Cartan equation).
Therefore, $\bar\6$ is a holomorphic structure operator.
Clearly, the holomorphic bundle $(B_1\oplus B_2, \bar\6)$ 
is isomorphic to $B$.

\hfill

Given an arbitrary harmonic form $\nu\in \c H^{0,1}(\End(B_{gr}))$,
the Maurer-Cartan equation \eqref{_MC_Equation_} 
will not, generally speaking, hold. The main purpose of
deformation theory is to find its solutions in terms
of the cohomology classes.

However, over a torus we may use the flatness
of harmonic representatives to obtain the solutions
of \eqref{_MC_Equation_}  in a straightforward way.
The following theorem is the main result of this Subsection.

\hfill

\theorem\label{_exists_nu_flat_Theorem_}
Let $T$ be a compact complex torus,
$B_1, ... B_n$ flat holomorphic Hermitian vector bundles
and $B$ a holomorphic vector bundle with a filtration
\[ 0 = E_0\subset E_1 \subset ...\subset E_n= B,\]
such that $E_i / E_{i-1} \cong B_i$. 
Then the following assertions are true.
\begin{description}
\item[(i)]
The holomorphic
structure on $B$ can be obtained as follows. 
Identify $B$ with $B_{gr} := \oplus B_i$
as $C^\infty$-bundle. Then there is a 
cohomology class $\nu$
\begin{equation}
\nu \in \bigoplus_{i>j} \Ext^1(B_i, B_j) \subset \Ext^1(B_{gr}, B_{gr})
\end{equation}
such that the holomorphic structure operator in $B$
is written as
\begin{equation}\label{_holo_via_nu_Equation_}
\bar \6 = \bar\6_{gr} + \nu_0, 
\end{equation}
where
$\bar\6_{gr}$ is the holomorphic structure
operator on $B_{gr}$, and $\nu_0 \in \Lambda^{0,1}(\End(B_{gr}))$
denotes the harmonic representative of $\nu$.

\item[(ii)] The class $\nu \in \oplus_{i>j} \Ext^1(B_i, B_j)$
is unique up to an automorphism of the bundle
$B_{gr} = \oplus B_i$ preserving the filtration
\[ 
  B_1 \subset B_1 \oplus B_2 \subset ... \subset \oplus_{i\leq k} B_i
  \subset ...
\]
\item[(iii)] Given an arbitrary cohomology class
\[ \nu \in \bigoplus_{i>j} \Ext^1(B_i, B_j) 
\]
with $\nu^2 =0$, we can reconstruct the holomorphic 
bundle $B$ as follows: $B$ is 
identified with $B_{gr}$ as a $C^\infty$-bundle,
and the holomorphic structure operator in $B$
is defined as in \eqref{_holo_via_nu_Equation_}.
\end{description}

\hfill

\remark
By \ref{_filtra_w_line_bu_Proposition_}, on a
complex compact  torus $T$ of complex dimension $d>2$ with
\[ H^{p,p}(T) \cap H^{2p}(T, \Z) =0, \text{\ \ for\ \ } 0<p<3,
\]
every vector bundle can be obtained this way.
Moreover, in this situation, 
the bundles $B_i$ can be chosen 1-dimensional.

\hfill

{\bf Proof of \ref{_exists_nu_flat_Theorem_}:}
Write the holomorphic structure operator on $B$ as
$\bar\6 = \bar\6_{gr} +\tilde \nu$,
where $\tilde\nu$ is a (0,1)-form with values
in \[ \oplus_{i>j} \Lambda^{0,1}(T, \Hom(B_i, B_j)).\]
The form $\tilde\nu$ satisfies the Maurer-Cartan
equation \eqref{_MC_Equation_}. However, this
form is, generally speaking, not harmonic.
Every automorphism $g\in \End B_{gr}$
acts on $\tilde \nu$ as 
$\tilde\nu\arrow g(\tilde\nu) + \bar\6_{gr}(g)$
(this is the well-known gauge action).
To produce $\nu$ with the properties described
in \ref{_exists_nu_flat_Theorem_} (i), we
need to find a correct gauge transform.

Consider the group $(\C^*)^n$ acting on
$B_{gr}$ by diagonal automorphisms, in such a way that
the $i$-th component $\alpha_i$ of $(\C^*)^n$ acts trivially
on $B_j\subset B_{gr}$ for $i\neq j$, and
as a multiplication by $\alpha_i$ on $B_i$.

We shall write the action of $(\C^*)^n$
on $\Lambda^{0,1}(\End(B_{gr})$ as follows.
Let 
\begin{align*} \tilde \nu &\in \Lambda^{0,1}(\End(B_{gr}) = 
   \oplus_{i,j}\Lambda^{0,1}(B_i, B_j),  \\ 
   \tilde\nu &:=  \sum_{i, j}\tilde\nu_{ij}, \ \ \ 
   \tilde\nu_{ij}\in \Lambda^{0,1}(B_i, B_j)
\end{align*}
If $\alpha \in(\C^*)^n$, $\alpha = \prod_i \alpha_i$,
then 
\[ \alpha(\tilde \nu) = \sum_{i, j} \alpha_i \alpha_j^{-1}
   \tilde\nu_{ij}
\]
The group $(\C^*)^n$ acts in this fashion on the solutions of
Maurer-Cartan equation, and maps every solution to an equivalent
one. If $\alpha_i \gg \alpha_j$ for all $i>j$, then
$\alpha$ maps 
\[ \tilde\nu\in\oplus_{i>j} \Lambda^{0,1}(T, \Hom(B_i, B_j))
\]
to a form which is arbitrarily small.

Consider the local deformation space $\Def(B_{gr})$ for
$B_{gr}$, constructed in \cite{_Siu_Trautmann_}.
The above argument implies that every
neighbourhood of the point $[B_{gr}]\in \Def(B_{gr})$
contains a bundle which is isomorphic to $B$.

Let $E$ be a holomorphic vector bundle over a compact
K\"ahler manifold. The local deformation space $\Def(E)$
can be constructed explicitly in terms of Massey products 
as follows. 

One can define the Massey products as obstructions
to constructing a solution of the Maurer-Cartan equation
(see e.g. \cite{_Babenko_Taimanov_}, or
\cite{_May_}, \cite{_Retakh_} for a more 
classical approach). Locally, $\Def(E)$ is 
embedded to the vector space $\Ext^1(E, E)$,
and the image of this embedding is a germ of all
vectors $\theta\in \Ext^1(E, E)$, such that
$\theta\wedge\theta=0$ and all the higher Massey
products of $\theta$ with itself vanish.

Fix such a vector $\theta\in \Ext^1(E, E)$.
We construct the corresponding 
vector bundle $E_\theta\in \Def(E)$
using the Hodge theory as follows (see e.g. 
\cite{_Verbitsky:Hyperholo_bundles_}).

Let $\theta_0 \in \c H^{0,1}(\Hom(E,E))$
be the harmonic representative of $\theta$. Using
induction, we define
\[ \theta_n := - \frac{1}{2} G_{\bar\6}\sum_{i+j=n-1} \theta_i \wedge \theta_j
\]
where $G_{\bar\6}$ is the Green operator inverting the
holomorphic structure operator
\[ \bar\6:\; \Lambda^{0,k}\otimes E\arrow  
    \Lambda^{0,1}\otimes E.
\]
on its image.
The Green operator $G_{\bar\6}$ is compact.
This can be used to show that 
for $\theta$ sufficiently small, the
series $\tilde \theta:= \sum \theta_i$
converges. The vanishing of Massey products
is equivalent to the following condition
\[ \bar\6 \theta_n = 
    - \frac{1}{2} \sum_{i+j=n-1} \theta_i \wedge \theta_j,
\]
which is apparent from the definition given in
\cite{_Babenko_Taimanov_}. In this case, we have
\[ \bar\6 \tilde \theta = 
   - \frac{1}{2} \tilde\theta\wedge\tilde\theta
\]
and $\tilde\theta$ is a solution of the Maurer-Cartan equation
\eqref{_MC_Equation_}.
Therefore, the operator $\bar\6_\theta= \bar\6+\tilde\theta$
satisfies $\bar\6_\theta^2=0$, and by Newlander-Nirenberg
theorem (\ref{_Newle_Nie_for_bu_Theorem_})
this operator defines a holomorphic structure
on $E$. On the deformation 
space $\Def(E) \subset \Ext^1(E, E)$,
the point $\theta$ corresponds to a bundle
$(E, \bar\6_\theta)$.

Now we return to holomorphic bundles over a compact torus
and the proof of \ref{_exists_nu_flat_Theorem_}. We obtain that
$B$ is given by some $\nu\in \Ext^1(B_{gr}, B_{gr})$. Since
$B$ and $\oplus B_{i}$ have compatible filtrations,
by functoriality we may assume that
\begin{equation}\label{_nu_nilpo_Equation_}
\nu \in \oplus_{i>j} \Ext^1(B_i, B_j).
\end{equation}
The higher Massey operations in $\Ext^*(B_{gr}, B_{br})$
vanish, because the bundle $B_{gr}$ is flat
(the same proof works as was used in \cite{_DGMS_};
see also \cite{_Goldman_Millson_}). Therefore,
$B$ can be reconstructed from $\nu$ 
for any $\nu$ such that the cohomology class 
$\nu\wedge\nu$ vanishes. Pick a harmonic representative
$\nu_0$ of $\nu$. Since $B_{gr}$ is flat, $\nu_0$ is
parallel. Therefore $\nu_0\wedge\nu_0$ is also parallel,
hence harmonic. We obtain that the cohomology
class of $\nu_0\wedge\nu_0$ vanishes if and only
if this form vanishes identically.

Starting from a bundle $B$ with a filtration
satisfying the assumptions of \ref{_exists_nu_flat_Theorem_},
we have constructed a cohomology class 
$\nu \in \oplus_{i>j} \Ext^1(B_i, B_j)$, with
$\nu\wedge\nu=0$. The harmonic representative $\nu_0$ of $\nu$
satisfies $\nu_0\wedge\nu_0=0$. Therefore, $\nu_0$ 
is a solution of Maurer-Cartan equation, and 
$\bar\6_{gr} +\nu_0$ is equivalent to 
the holomorphic structure operator of $B$.
This proves \ref{_exists_nu_flat_Theorem_} (i).
\ref{_exists_nu_flat_Theorem_} (ii) is clear 
by functoriality of our construction, and
\ref{_exists_nu_flat_Theorem_} (iii) is obvious.
\endproof

\hfill

We also obtained the following corollary.

\hfill

\corollary\label{_bun_admits_flat_Corollary_}
Let $B$ be a holomorphic vector bundle on a compact
torus $T$, $\dim_\C T>2$,
\[ 
H^{p,p}(T) \cap H^{2p}(T, \Z) =0, \text{\ \ for\ \ } 0<p< 3
\]
Then $B$ admits a flat connection compatible with the holomorphic structure.

\hfill

{\bf Proof:} By \ref{_filtra_w_line_bu_Proposition_},
$B$ admits a filtration satisfying the assumptions of 
\ref{_exists_nu_flat_Theorem_}. Let $\nu_0\in \Lambda^{0,1}(\End(B_{gr}))$
be the $(0,1)$-form which defines the holomorphic structure
on $B$ as in \eqref{_holo_via_nu_Equation_}:
\[ \bar \6 = \bar\6_{gr} + \nu_0, \]
and let $\nabla_{gr}$ be the standard 
flat Hermitian connection on $B_{gr}$.
Then the connection $\nabla:= \nabla_{gr}+ \nu_0$ is flat.
By construction, the $(0,1)$-part of $\nabla$
is equal to the holomorphic structure operator on $B$.
\endproof

\hfill

The form $\nu_0$ is called {\bf the Higgs field of the holomorphic
bundle $B$}.


\section{Yoneda extension in abelian categories}
\label{_Yoneda_Section_}


In this Section we study abelian categories of finite length,
in quite an abstract setting. Further on, these results
are applied to the category of reflexive sheaves on a hyperk\"ahler
manifold.

\subsection{Abelian categories of cohomological dimension $\leq 1$.}

\definition
Let $\c C$ be an abelian category.
An object $B\in \c C$ is called {\bf simple}
if it has no proper sub-objects:
for all $B'\subset B$, either $B'=B$
or $B'=0$. An object $B$ is called {\bf semisimple}
if $B$ is a direct sum of simple objects.

An abelian category $\c C$ is called {\bf of finite length } if 
every object $B\in \c C$ is a finite extension
of simple objects. Equivalently, $\c C$ 
is of finite length if any increasing or decreasing chain of 
sub-objects of $B$ stabilizes, for all 
$B\in \c C$.\footnote{One also says {\bf $\c C$ 
satisfies the ascending and descending chain condition.}}

\hfill

Let $\c C$ be an abelian category. Denote by $\Ext^1(B, B')$
the group of Yoneda extensions from $B$ to $B'$.
Given an exact sequence
\[ 0 \arrow B_1 \arrow B_2 \arrow B_3 \arrow 0
\]
we have the exact sequences
\begin{equation}\label{_Ext_long_exa_seq2_Equation_}
\begin{aligned}
0 &\arrow \Hom(B_3, B) \arrow \Hom(B_2, B) \arrow
\Hom(B_1, B) \arrow\\ {}& \arrow \Ext^1(B_3, B) \arrow \Ext^1(B_2, B) \arrow
\Ext^1(B_3, B) 
\end{aligned}
\end{equation}
and
\begin{equation}\label{_Ext_long_exa_seq1_Equation_}
\begin{aligned}
0 &\arrow \Hom(B, B_1) \arrow \Hom(B, B_2) \arrow
\Hom(B, B_3)\arrow \\ {}& \arrow \Ext^1(B, B_1) \arrow \Ext^1(B, B_2) \arrow
\Ext^1(B, B_3).
\end{aligned}
\end{equation}

The following elementary claim is well known

\hfill

\claim\label{_coho_1_Claim_}
Let $\c C$ be an abelian category.
Then the following conditions are equivalent:
\begin{description}
\item[(i)] The exact sequences \eqref{_Ext_long_exa_seq2_Equation_},
\eqref{_Ext_long_exa_seq1_Equation_} are left exact
\item[(ii)] The second Yoneda extension group vanishes identically
\item[(iii)] Any exact sequence of length 4 splits onto a direct sum
of two exact sequences of length 3.
\end{description}
\endproof

\hfill

\definition
Let $\c C$ be an abelian category satisfying 
the conditions (i)-(iii) of \ref{_coho_1_Claim_}.
We say that $\c C$ has cohomological dimension $\leq 1$.

\hfill

The following lemma is quite easy.

\hfill

\lemma \label{_fini_length_coho_1_Lemma_}
Let $\c C$ be an abelian category of finite length.
Assume that the exact sequences \eqref{_Ext_long_exa_seq2_Equation_},
\eqref{_Ext_long_exa_seq1_Equation_} are left exact
for all semisimple $B, B_3$.
Then $\c C$ is of cohomological dimension $\leq 1$. 

\hfill

{\bf Proof:} Define the {\bf length} of an object $F\in \C$ 
as a length of a minimal filtration on $B$ with semisimple
associated graded factors. We denote the length of $F$ by $l(F)$.

We prove \ref{_fini_length_coho_1_Lemma_} by induction
by $l(B)$, $l(B_i)$. By assumptions of 
\ref{_fini_length_coho_1_Lemma_},
the second Yoneda extension
$\Ext^2(B, B')$ vanishes 
for all semisimple $B, B'$.
We need to show that $\Ext^2(B, B')=0$
for arbitrary $B, B'$. Using induction, 
we may assume that this group vanishes 
for all $B, B'$ with $l(B)<n, l(B')<m$.
Write an exact sequence
\[ 0 \arrow B'' \arrow B' \arrow B''' \arrow 0
\]
with $l(B''), l(B''') < l(B')$. From the long exact sequence
\[
...\arrow \Ext^2(B,B'') \arrow \Ext^2(B,B')\arrow \Ext^2(B,B''') \arrow ... 
\]
and the induction assumption we obtain that $\Ext^2(B,B')=0$.
\endproof

\hfill

\lemma\label{_cohomo_dime_subcate_Lemma_}
Let $\c C \subset \c C_0$ be a full abelian subcategory of an
abelian category $\c C_0$. Assume that $\c C_0$ is of cohomological
dimension $\leq 2$. Then $\c C$ is of cohomological dimension 
$\leq 1$.

\hfill

{\bf Proof:} Take a 4-term exact sequence in $\c C$
\begin{equation}\label{_4-term_exa_Equation_}
0 \arrow F_1\arrow F_2 \arrow F_3 \arrow F_4\arrow 0.
\end{equation}
Since $\c C_0$ has cohomological dimension $\leq 1$,
\eqref{_4-term_exa_Equation_} splits in $\c C_0$ onto
a direct sum of two 3-term exact sequences. Since
$\c C$ is full in $\c C_0$, all direct summands
of objects in $\c C$ also belong in $\c C$. Therefore,
these 3-term exact sequences belong in $\c C$, and
\eqref{_4-term_exa_Equation_} splits in $\c C$ onto
a direct sum of two 3-term exact sequences. \endproof

\hfill

\subsection{Functors on abelian categories of finite length}

The main result of this Section is the following 
elementary proposition

\hfill

\proposition
\label{_functo_inducing_isomo_on_semisi_Proposition_}
Let $\c C \stackrel \gamma \arrow \c C'$ be a functor of
abelian categories of finite length 
satisfying the following properties
\begin{description}
\item[(i)] $\gamma$ induces an equivalence on the categories
of semisimple objects in $\c C$, $\c C'$
\item[(ii)] $\gamma$ induces an isomorphism on the Yoneda
$\Ext^1$-groups
\[ \Ext^1(B_1, B_2) \arrow \Ext^1(\gamma(B_1), \gamma(B_2)),
\]
where $B_1$, $B_2$ are semisimple objects in $\c C$.
\item[(iii)] The categories $\c C$, $\c C'$ are of cohomological
  dimension $\leq 1$.
\end{description}
Then $\gamma$ is an equivalence.

\hfill

{\bf Proof:}\footnote{This
proof was suggested by D. Kaledin.} 
The proof is obtained via the trivial diagram-chasing
argument. As in the proof of \ref{_fini_length_coho_1_Lemma_},  
we use an induction by the length of $B$. Let $B, B'\in \c C$
be objects of length $n$. Consider an exact sequence
\[ 0 \arrow B_0 \arrow B \arrow B_1 \arrow 0
\]
with $l(B_1)= n-1$ and $B_0$ semisimple. There is a
long exact sequence
\[
0 \to \Hom(B', B_0) \to \Hom(B', B) \to 
\Hom(B', B_1) \to \Ext^1(B', B_0) \to ...
\]
Applying $\gamma$, we obtain a commutative diagram with
exact rows

{\tiny
\minCDarrowwidth1.6pc
\begin{equation}\label{_long_exa_Phi_Ext^1_Equation_}
\begin{CD}
0 @>>> \Hom(B', B_0) @>>> \Hom(B', B) @>>> 
\Hom(B', B_1)\arrow \\
&& @V{\gamma} VV @V{\gamma}VV  @V{\gamma}VV \\
0 @>>> \Hom(\gamma(B'), \gamma(B_0)) @>>> \Hom(\gamma(B'), \gamma(B)) @>>> 
\Hom(\gamma(B'), \gamma(B_1)) \arrow \\[4mm]
&&  \arrow \Ext^1(B', B_0) @>>> \Ext^1(B', B) @>>> \Ext^1(B', B_1) \arrow \\
&& @V{\gamma} VV @V{\gamma}VV @V{\gamma}VV \\
&& \arrow \Ext^1(\gamma(B'), \gamma(B_0))
@>>> \Ext^1(\gamma(B'), \gamma(B))@>>>\Ext^1(\gamma(B'), \gamma(B_1)) \arrow  \\[4mm]
&&&& &&  \arrow  \Ext^2(B', B_0)\\
&&&&&& @V{\gamma}VV  \\
&&&&&&  \arrow  \Ext^2(\gamma(B'), \gamma(B_0))
\end{CD}
\end{equation}
}
Since $\c C$, $\c C'$ is of cohomological dimension $\leq 1$,
we have $\Ext^2(B', B_0)=\Ext^2(\gamma(B'), \gamma(B_0))=0$.
Using induction, we may assume that all vertical arrows of
\eqref{_long_exa_Phi_Ext^1_Equation_} are isomorphisms,
except, possibly, the maps 
\[ \Ext^1(B', B)\stackrel \gamma\arrow \Ext^1(\gamma(B'), \gamma(B))\]
and \[ \Hom(B', B)\stackrel \gamma\arrow \Hom(\gamma(B'), \gamma(B)).\]
The snake lemma implies immediately that 
these two maps are isomorphisms as well.
We proved that $\gamma$ is full and faithful
and induces an isomorphism on $\Ext^1$-groups.
It remains to show that any object of $\c C'$
is isomorphic to $\gamma(B)$ for some $B\in \c C$.

Consider an exact sequence
\[ 0 \arrow B_0 \arrow B \arrow B_1 \arrow 0
\]
with $l(B_1)= n-1$ and $B_0$ semisimple. 
Using induction, we may assume that every object of length
$\leq n-1$ is isomorphic to an object from the image of $\gamma$. 
Then, $B_0$ and $B_1$ belong to the image of $\gamma$.
Since $B$ is represented by $\Ext^1(B_1, B_0)$, 
and $\gamma$ induces isomorphism on $\Ext^1$-groups,
$B$ also belongs to the image of $\gamma$.
This proves \ref{_functo_inducing_isomo_on_semisi_Proposition_}.
\endproof

\hfill

The same arguments as prove
\ref{_functo_inducing_isomo_on_semisi_Proposition_}
can also be used to prove the following corollary

\hfill

\corollary \label{_embe_on_Ext_1_hence_full_faith_Corollary_}
Let $\c C \stackrel \gamma \arrow \c C'$ be a functor of
abelian categories of finite length. Assume that
$\gamma$ induces isomorphism on the respective
subcategories of semisimple objects. Assume,
moreover, that $\gamma$ induces
a monomorphism on the Yoneda $\Ext^1$-groups
between semisimple objects. Then $\gamma$ is full and failthful.

\endproof


\section{Reflexive sheaves on twistor spaces}
\label{_refle_equi_Section_}


The aim of this section is the following theorem.

\hfill

\theorem \label{_refle_shea_isomo_Theorem_}
Let $M$ be a hyperk\"ahler K3 surface, or a hyperk\"ahler compact
torus. Assume that $M$ is generic in the sense of having
no non-trivial $SU(2)$-invariant integer cycles:
\[ (H^p(M)_{SU(2)-inv})\cap H^p(M, \Z) =0, \text{ \ \ \ for all\ \ \ }
   0<p< \dim_\R M.
\]
Consider generic\footnote{In the sense of \ref{_generic_manifolds_Definition_}} 
complex structures $L, L'$ induced by the hyperk\"ahler structure,
and let $\c C^r(M,L)$, $\c C^r(M,L')$ be the categories of 
reflexive sheaves on $(M,L)$ and $(M, L')$. Then the
categories $\c C^r(M,L)$, $\c C^r(M,L')$ are equivalent.
Moreover, this equivalence can be chosen in a canonical way
if we chose another induced complex structure $I$,
$I\neq L, L'$.

\hfill

\ref{_refle_shea_isomo_Theorem_}
is proven in Subsection \ref{_R_L_equi_refle_for_surfa_Subsection_} 
for $M$ a surface,
and in Subsection \ref{_C^tw_for_torus_Subsection_} 
for $M$ a compact torus,
$\dim_\C M >2$.

\subsection{Category $\c C^{tw}_I$: the definition}
\label{_C^tw_I_defini_Subsection_}

We work in assumptions of \ref{_refle_shea_isomo_Theorem_}.
Consider $(M,L)$ as a submanifold in $\Tw(M)$,
$(M, L) = \pi^{-1}(L)$. 
Consider the category 
\[
   \c C^r:= \lim_\arrow \c C^b(\Tw(M \backslash S))
\]
where $\c C^b(\Tw(M \backslash S))$ is the category of vector bundles
on $\Tw(M \backslash S)$ which can be extended 
to $\Tw(M)$, and the limit is taken over all increasing
sequences of finite subsets $S \subset M$. By \ref{_refle_pushfor_Lemma_}, 
we may think of the objects of $\c C^r$ as of reflexive sheaves
on $\Tw(M)$; the space of morpfisms of objects of $\c C^r$
is isomorphic to the space of morphisms of the corresponding
reflexive sheaves. 

We are going to embed the category
of reflexive sheaves on $(M,L)$ into the category 
$\c C^r$ in such a 
way that the image is independent from the choice of $L$.
Varying $L$, we obtain identifications between
$\c C^r(M,L)$ for various generic induced complex 
structures $L$. From the definition of $\c C^r$
it follows that the natural restriction functor
$R_L(\cdot) := (\cdot) \restrict{(M,L)}$
maps $\c C^r$ to the category
$\c C^r(M,L)$ of reflexive sheaves
on $(M,L)$. 

\hfill

Fix an induced complex structure $I\neq L$ on $M$.
We define the subcategory
$\c C^{tw}_I \subset \c C^r$,
in such a way that the restriction functor
\[ R_L:\; \c C^{tw}_I \arrow \c C^r(M,L)\]
is an equivalence, for all generic induced
complex structures $L\neq I$. 

\hfill

In \cite{_Li_Yau_}, a theory of stable sheaves was constructed
for some non-K\"ahler manifolds. In \cite{_NHYM_}, 
we proved that this theory can be applied to the twistor
spaces. In particular, one may speak of
semistable sheaves on $\Tw(M)$ and of 
Jordan-Holder filtration, that
is, a filtration with stable associated
graded factors. As it happens in the usual case,
the associated graded sheaf of the Jordan-Holder filtration
is independent from the choice of such filtration.

\hfill

\definition\label{_C_I^tw_Definition_}
Let $M$ be a compact hyperk\"ahler torus or a K3 surface,
which is generic in the sense of having
no non-trivial $SU(2)$-invariant integer cycles:
\[ (H^p(M)_{SU(2)-inv})\cap H^p(M, \Z) =0, \text{ \ \ \ for all\ \ \ }
   0<p< \dim_\R M.
\]
The category
$\c C^{tw}_I$ is a full subcategory of the category
of $\c C^r$ reflexive sheaves on $\Tw(M)$, defined as follows. 
Consider a reflexive sheaf $F\in \c C^r$.
Then $F$ belongs to $\c C^{tw}_I$ iff the following
conditions hold
\begin{description}
\item[(i)] The sheaf $F$ is semistable. The restriction of $F$
to $(M, I)\subset \Tw(M)$ has a polystable reflexization:
\[ \bigg(F\restrict{(M,I)}\bigg)^{**} \cong \oplus B_i
\]
where $B_i$ are hyperholomorphic bundles on $(M,I)$

\item[(ii)] 
There is a finite set $A\subset M$ such that
$\c F\restrict{\Tw(M) \backslash (A \times \{I\})}$ is a bundle
(that is, all singularities of $F$ are sitting in the finite
set $A \times \{I\}\subset (M, I) \subset \Tw(M)$).

\item[(iii)] Consider a Jordan-Holder filtration of $F$
\[  
   0 = F_0\subset F_1\subset .... \subset F_n = F
\]
with all $F_i$ reflexive, and the quotient sheaves $F_i/F_{i-1}$
stable. Then the quotient sheaves $F_i/F_{i-1}$
are non-singular on $\Tw(M \backslash A)$.
For all $i$, the reflexization of $F_i / F_{i-1}$
is isomorphic to $\Tw(B_i)$, where $\Tw(B_i)$
is the twistor transform of the bundle $B_i$
considered in (i).

\item[(iv)] Let $x\in A$, and let $\Tw(M)_{x, I}$ be the germ
of a neighbourhoof of $s_x \backslash I$ within
$\Tw(M) \backslash (M,I)$ (see 
subsection \ref{_twi_spa_local_Subsection_}).
Then the sheaf $F\restrict{\Tw(M)_{x, I}}$ can be trivialized
over $s_x \backslash I$ by the map
\[
\tau:\; F\restrict {\Tw(M)_{x_0, I}}\arrow \xi^*F_S.
\]
(see \eqref{_E_loca_triv_Equation_} for details
and notation). Moreover, this trivialization induces
a trivialization on all subsheaves $F'\subset F$.
Whenether the quotient of two subsheaves
$F_1/F_2$ has a reflexization isomorphic to
$\Tw(B)$, then $\tau$ induces the same trivialization 
as \eqref{_E_loca_triv_Equation_} gives.

\end{description}

We consider objects of $\c C^{tw}_I$ as reflexive sheaves on
$\Tw(M)$.

\hfill

\lemma
The category $\c C^{tw}_I$ is abelian.

\hfill

{\bf Proof:}
A full additive subcategory of an abelian category is abelian
if and only if for any morphism $\phi$ of this subcategory,
the kernel and the cokernel of $\phi$ belongs to this
subcategory. Let $F, F' \in \c C^{tw}_I$, and let
$\phi:\; F \arrow F'$ be any morphism.
Since $F$, $F'$ is semistable, 
the sheaves $\ker\phi$ and $\coker \phi$ are 
also semistable. Moreover, the Jordan-Holder filtration
on $F$, $F'$ induces the Jordan-Holder filtration on
$\ker\phi$ and $\coker \phi$, in such a way that the
associated graded sheaf of the Jordan-Holder filtration
on $\ker\phi$, $\coker \phi$ is a direct sum component 
of the associated graded sheaf of the Jordan-Holder filtration
on $F$, $F'$. This immediately implies that
the conditions (i)-(iii) of \ref{_C_I^tw_Definition_}
hold for $\ker\phi$ and $\coker \phi$.
The condition (iv) holds for
$\ker\phi$ and $\coker \phi$
tautologically.
\endproof

\hfill

\claim\label{_exte_in_C^tw_I_in_terms_of_loca_Claim_}
Suppose that $\dim_\C M =2$.
Let $F\in \c C^r$ be a sheaf satisfying
the conditions (i)-(iii) of \ref{_C_I^tw_Definition_}.
Suppose that the Jordan-Holder filtration on $F$ has length 2,
and let $B_1$, $B_2$ be the corresponding hyperholomorphic
vector bundles. We have an exact sequence of vector bundles
on $\Tw(M\backslash A)$
\[ 0 \arrow \Tw(B_1)\restrict{\Tw(M\backslash A)} \arrow F\restrict{\Tw(M\backslash A)} \arrow \Tw(B_2)\restrict{\Tw(M\backslash A)}\arrow 0.
\]
Let 
\[ \nu \in \Ext^1\bigg(\Tw(B_1)\restrict{\Tw(M\backslash A)}, \Tw(B_2)\restrict{\Tw(M\backslash A)}\bigg)\cong H^1(\Tw((B_1)^* \otimes B_2))
\]
be the corresponding element of the extension group.
Denote by $B$ the bundle $(B_1)^* \otimes B_2$.
Consider the image of $\nu$
inside the group 
\begin{align*} 
\ker\bigg( 
\bigoplus_{x_i\in A} & \Hom_\C(\hat\calo_{s_{x_i}},\C)\otimes_{\calo_{\Tw(M)}}
  \Tw(B) \otimes_{\calo(\C P^1)} \calo(2) \\
  {}& \arrow R^2\pi_* \Tw(B) \bigg)
\end{align*}
induced by the exact sequence
\eqref{_cohomo_w_suppo_on_twi_long_exa_Equation}.
Then (iv) holds if and only if
$\nu$ is compatible with the trivialization
constructed in \eqref{_loca_coho_sheaf_trivialized_Equation_}.

\hfill

{\bf Proof:} Let $x\in A$ and let $U$ be a small neighbourhood
of \[ s_x \backslash I \subset \Tw(M) \backslash (M, I). \]
Assume that the fibers of the natural projection
$\pi:\; U \arrow \C P^1$ are Stein.
Consider the manifold $U^0:= U \backslash s_x$. The
bundle $F\restrict {U_0}$ is an extension of
$B_1\restrict {U^0}$ and $B_2\restrict {U^0}$.
This extension might be non-trivial; its non-triviality
is controlled by the group 
\[
P:= \ker\bigg( 
\Hom_\C(\hat\calo_{s_{x}},\C)\otimes_{\calo_{\Tw(M)}}
  \Tw(B) \otimes_{\calo(\C P^1)} \calo(2) \arrow R^2\pi_* \Tw(B) \bigg).
\]
For any $\nu \in P$,
and any $J\in s_x \backslash I$, consider the restriction
$F^J$ of $F$ to $\pi^{-1}(J) \subset U^0$. Then 
$F^J\restrict{U^0}$ is an extension of $B_1^J\restrict{U^0}$ and $B_2^J\restrict{U^0}$. Since $\pi^{-1}(J)$
is a 2-dimensional Stein manifold with a point deleted, this
extension might be non-trivial; it is controlled
by the restriction of $\nu$ to $J$, which
is a vector in
\[
\ker\bigg( 
\Hom_\C(\hat\calo_x,\C)\otimes_{\calo_{(M, J)}}
  B^J\otimes K_{(M,J)}
  \arrow H^2(B^J) \bigg).
\]
(see \ref{_Gro_support_Proposition_}).
The trivialization preserves $F$ and the subsheaf $F_1 \subset F$, 
hence it preserves  the cohomology class of the extension $\nu$.
Therefore, $\nu$ is trivialized.
The converse is also clear, because $F\restrict U$
is a pushforward of $F\restrict {U^0}$ (\ref{_normal_shea_Lemma_}). 
\endproof

\hfill

\remark
If $A$ is empty, then for all $F\in \c C^{tw}_I$,
the sheaf $F$  
and all the sheaves $F_i$ and $F_i/F_j$ of
\ref{_C_I^tw_Definition_} are non-singular.
In this case, $F$ admits a filtration by holomorphic bundles $F_i$
such that $F_i/F_{i-1}$ is isomorphic to $\Tw(B_i)$, 
where $B_i\cong (F_i/F_{i-1})\restrict{(M,I)}$, 
and $F\restrict{(M,I)} \cong \oplus B_i$. Moreover, every
bundle $F$ on $\Tw(M)$ admitting such a filtration 
belongs to $\c C^{tw}_I$. If $M$ is a torus,
$\dim_\C M >2$, then $A$ is always empty,
because all reflexive sheaves on $\Tw(M)$ 
are locally free, and the associated graded sheaves of the
Jordan-Holger filtration are also locally free
(this result can be proven in the same fashion as 
one proves \ref{_refle_she_on_torus_bundles_Lemma_}
and \ref{_filtra_w_line_bu_Proposition_}).

\hfill

\ref{_refle_shea_isomo_Theorem_} is immediately implied by
the following theorem.

\hfill

\theorem\label{_C^tw_I_isomo_Theorem_}
Let  $M$ be a compact hyperk\"ahler torus or a K3 surface,
$I$ an induced complex structure, and $L\neq I$ a
generic induced complex structure.
Assume that $M$ is generic in the sense of having
no non-trivial $SU(2)$-invariant integer cycles:
\[ (H^p(M)_{SU(2)-inv})\cap H^p(M, \Z) =0, \text{ \ \ \ for all\ \ \ }
   0<p< \dim_\R M.
\]
Consider the category $\c C^{tw}_I$ defined above,
and let 
\begin{equation}\label{_R_L_restri_Equation_}
R_L:\; \c C^{tw}_I \arrow \c C^r(M,L)
\end{equation}
be the restriction map from $\c C^{tw}_I$
to the category $\c C^r(M,L)$ of reflexive sheaves
on $(M,I)$. Then \eqref{_R_L_restri_Equation_} is an equivalence
of categories.

\hfill

\ref{_C^tw_I_isomo_Theorem_}
is proven in Subsection \ref{_R_L_equi_refle_for_surfa_Subsection_} 
for $M$ a surface,
and in Subsection \ref{_C^tw_for_torus_Subsection_} 
for $M$ a compact torus of dimension
greater than 2.

\hfill

The definition of $\c C^{tw}_I$ is motivated by the
following heuristic consideration. 
We need to have a subcategory of the category
$\c C^r$
of reflexive sheaves on $\Tw(M)$ which is isomorphic
to $\c C^r(M,L)$. The category of semisimple objects of
$\c C^r(M,L)$ is embedded to the category of stable bundles
on $\Tw(M)$ in a natural way via the twistor transform
(see Subsection \ref{_twi_fo_for_bu_Subsection_}).
By \ref{_direct_ima_twi_Proposition_}, for any
pair of hyperholomorphic bundles $B, B'$ on $(M,L)$, 
we have
\[ 
\Ext^1(\Tw(B), \Tw(B')) = \C^2\otimes_\C \Ext^1(B, B').
\]
To obtain within  $\c C^r$ a subcategory isomorphic
to $\c C^r(M,L)$, we need to ``kill off'' half of the
extensions. This is performed quite easily by requiring that
any extension from $\c C^{tw}_I$ splits over $(M,I)$.
Indeed, by  \ref{_direct_ima_twi_Proposition_},
the extensions within the category of coherent sheaves
on $\Tw(M)$ are
just sections of $\calo(1) \otimes \Ext^1(B, B')$.
The sections of 
$\calo(1) \otimes \Ext^1(B, B')$ which vanish at $I$
are in 1-to-1 correspondence with $\Ext^1(B, B')$.

If we deal with vector bundles, the 
definition of $\c C^{tw}_I$ can be reformulated
in a more concise fashion, by saying that
$F$ belongs to $\c C^{tw}_I$ if $F$ is an extension
of bundles of type $\Tw(B_i)$ which splits at $(M,I)$.

The rest of the definition of $\c C^{tw}_I$ 
adapts the same type of reasoning by accomodating
singularities into the picture.

\subsection{Reflexive sheaves on generic compact torus and
the category $\c C^{tw}_I$}
\label{_C^tw_for_torus_Subsection_}

Let $T$ be a hyperk\"ahler compact torus, $\dim_\C T >2$. 
Assume that
all $SU(2)$-invariant integer classes
$\eta \in H^p_{SU(2)-inv}(T) \cap H^p(T, \Z)$,
$0 < p <\dim_\R T$ vanish. We are going to prove
\ref{_C^tw_I_isomo_Theorem_} for the manifold $T$.

All reflexive coherent sheaves on $T$ 
are bundles (\ref{_cohe_on_toru_Theorem_}). 
The same argument applied to the twistor space
shows that all $F\in \c C^{tw}_I$ are also bundles. 
By \ref{_filtra_w_line_bu_Proposition_},
all bundles on $T$ admit a filtration
by vector bundles, with associate graded sheaves polystable
and locally free. The same argument implies that
all $E\in \c C^{tw}_I$ admit a holomorphic filtration
\[ 0 = E_0\subset E_1 \subset ...\subset E_n= E,\]
with $E_i/ E_{i-1}= \Tw(B_i)$, where $B_i$ is a
hyperholomorphic bundle on $(T,I)$, such that
$E\restrict{(T,I)}\cong \oplus B_i$. Conversely,
every such $E$ belongs to $\c C_I^{tw}$.

Let $L$ be a generic induced complex structure on $T$.
Denote by $\c C^b(T,L)$ the category
of holomorphic vector bundles on $(T,L)$.
By \ref{_cohe_on_toru_Theorem_}, all reflexive sheaves
on $(T,L)$ are vector bundles. Clearly, the
category $\c C^b(T,L)= \c C^r(T,L)$ is abelian. 
Consider the restriction functor
\begin{equation}\label{_R_L_on_T_Equation_}
R_L:\; \c C^{tw}_I \arrow \c C^b(T,L).
\end{equation}
To prove \ref{_C^tw_I_isomo_Theorem_} for $T$,
we need to show that \eqref{_R_L_on_T_Equation_}
is an equivalence. 

As we have indicated above, all simple objects of $\c C^{tw}_I$
are of form $\Tw(B)$, where $B$ is a stable bundle on $T$. The
simple objects of $\c C^b(T,L)$ are stable bundles as well. 
Therefore, \eqref{_R_L_on_T_Equation_}
induces an isomorphism of respective 
subcategories of semisimple objects. 
We arrive in the situation similar
to that described in
\ref{_embe_on_Ext_1_hence_full_faith_Corollary_}.
To apply \ref{_embe_on_Ext_1_hence_full_faith_Corollary_},
we need to compare the Yoneda extensions between the
semisimple objects in $\c C^{tw}_I$ and $\c C^b(T, L)$.

Let $F, F'$ be bundles on $\Tw(T)$,
$F = \Tw(B)$, $F'= \Tw(B')$, with $B$, $B'$
hyperholomorphic bundles on $T$. Denote by $B_L$, $B_L'$ the bundles
$B$ and $B'$ considered as holomorphic vector bundles
on $(T, L)$:
\[ B_L= R_L(F), \ \ \  B_L' = R_L(F').
\]
By \ref{_direct_ima_twi_Proposition_}, we have
\begin{equation}\label{_Ext^1_O(1)_Equation_}
\Ext^1(F, F') = \Ext^1(B_L, B'_L) \otimes_\C \c H^0(\C P^1, \calo(1)).
\end{equation}
An extension 
\[ 
0 \arrow F' \arrow F'' \arrow F \arrow 0
\]
belongs to $\c C^{tw}_I$ if and only if it splits on $(T,I) \subset \Tw(T)$.
This allows us to compute the group 
$\Ext^1_{\c C^{tw}_I}(F, F')$
of Yoneda extensions in $\c C^{tw}_I$ from $F$ to $F'$.
Looking at \eqref{_Ext^1_O(1)_Equation_},
we obtain that
\[ 
\Ext^1_{\c C^{tw}_I}(F, F') = 
\Ext^1(B_L, B'_L) \otimes_\C \c H^0(\C P^1, \calo(1))_I,
\]
where $\c H^0(\C P^1, \calo(1))_I$
denotes the space of all sections of $\calo(1)$
vanishing at $I$. However, the group
$\c H^0(\C P^1, \calo(1))_I$ is clearly 1-dimensional,
and this gives
\[ 
\Ext^1_{\c C^{tw}_I}(F, F') \cong
\Ext^1(B_L, B'_L).
\]
The functor \[ R_L:\; \c C^{tw}_I \arrow \c C^b(T,L)\] evaluates 
a section $\eta\in \c H^0(\C P^1, \calo(1))_I$ at $L\in \C P^1$, hence
it induces an isomorphism on the first Yoneda extensions
$\Ext^1$, assuming that $L\neq I$.
Applying \ref{_embe_on_Ext_1_hence_full_faith_Corollary_}, 
we obtain that  \[ R_L:\; \c C^{tw}_I \arrow \c C^b(T,L)\]
is full and faithful.

\hfill

To finish the proof of \ref{_C^tw_I_isomo_Theorem_} for
$T$, it remains to show that $R_L$ is surjective on the 
set of equivalence classes of objects. Let $B$ be a holomorphic 
bundle on $(T,L)$. We need to show that there is a bundle
$F \in \c C^{tw}_I$ such that $B \cong F\restrict{(M,L)}$.
Using \ref{_exists_nu_flat_Theorem_}, we write 
\[
B = (B_{gr}, \bar\6_{gr}+\nu),
\]
where $B_{gr}$ is a flat Hermitian bundle on $T$,
$\bar\6_{gr}$ the holomorphic structure operator
on $B_{gr}$ and $\nu\in \Lambda^{0,1}(T, \End(B_{gr}))$ 
the Higgs field. 
Denote by $\nabla_{gr}$ the flat 
Hermitian connection given on $B_{gr}$.
Consider the flat connection
\[ \nabla_\nu := \nabla_{gr}+\nu
\]
on $B$ (\ref{_bun_admits_flat_Corollary_}).
Lifting $(B, \nabla_\nu)$ as in \ref{_twi_functor_Theorem_}, 
we obtain a holomorphic bundle $F_\nu$ on $\Tw(T)$.
Clearly, $F_\nu \in \c C^{tw}_{-L}$, and $F_\nu\restrict{(M,L)}=B$.
This proves \ref{_C^tw_I_isomo_Theorem_} for $M=T$ and
$L=-I$. To prove \ref{_C^tw_I_isomo_Theorem_} for arbitrary
$L$, we write another flat connection on $B$, as follows.

Let $\alpha \in \C$ be an arbitrary number. Consider the
1-form $\nu_\alpha := \nu + \alpha \bar\nu$.

\hfill

\claim\label{_nu_alpha_MC_Claim_}
In the above assumptions, the form $\nu_\alpha = \nu + \alpha \bar\nu$
satisfies the Maurer-Cartan equation:
\[
\nabla_{gr}(\nu_\alpha) + 2 \nu_\alpha\wedge\nu_\alpha=0
\]

\hfill

{\bf Proof:}
Clearly, $\nu_\alpha$ is parallel with respect
to the connection $\nabla_{gr}$. To show that
$\nu_\alpha\wedge\nu_\alpha=0$, we need only to prove
that $\nu \wedge\bar\nu + \bar\nu \wedge\nu =0$. 
However,
\begin{align*} (\nu+\bar\nu)\wedge(\nu+\bar\nu) & 
   =  \nu\wedge\nu+ \bar\nu\wedge\bar\nu
   + \nu \wedge\bar\nu + \bar\nu \wedge\nu \\ 
   {} &= \nu \wedge\bar\nu + \bar\nu \wedge\nu,
\end{align*}
hence to show $\nu_\alpha\wedge\nu_\alpha=0$
we need only to prove that 
$(\nu+\bar\nu)\wedge(\nu+\bar\nu)=0$.
Write the connection
$\nabla_{r}:= \nabla_{gr} + \nu+\bar\nu$.
This connection is clearly  compatible 
with the Hermitian structure on $B_{gr}$.
The $(0,1)$-part of $\nabla_r$ is equal to
$\bar\6_{gr} + \nu$, hence its square vanishes,
and $\nabla_r$ is a holomorphic Hermitian connection.
The curvature of $\nabla_r$ is written as
 \begin{align*}
\nabla_r^2 = \nabla_r(\nu+\bar\nu) +(\nu+\bar\nu)\wedge(\nu+\bar\nu) & 
   =  \nu\wedge\nu+ \bar\nu\wedge\bar\nu
   + \nu \wedge\bar\nu + \bar\nu \wedge\nu \\ 
   {} &= \nu \wedge\bar\nu + \bar\nu \wedge\nu,
\end{align*}
because $\nu+\bar\nu$ is parallel. Therefore, to prove
$\nu \wedge\bar\nu + \bar\nu \wedge\nu =0$ it means to prove
that $\nabla_r$ is flat. However, the form
$\nu \wedge\bar\nu + \bar\nu \wedge\nu$
is  parallel, hence harmonic. By \ref{_harmonic_curv_Remark_},
the connection $\nabla_r$ is Yang-Mills.
Since the Chern classes $c_1(B_{gr})$,
$c_2(B_{gr})$ vanish, any Yang-Mills connection on $B$
is flat (\cite{_Lubcke_}, \cite{_Simpson_}).
Therefore, $\nabla_r^2=0$, and we have
$\nu \wedge\bar\nu + \bar\nu \wedge\nu =0$.
We proved \ref{_nu_alpha_MC_Claim_}.
\endproof

\hfill

By \ref{_nu_alpha_MC_Claim_}, 
the connection $\nabla_{\alpha}:= \nabla_{gr}+ \nu + \alpha \bar\nu$
is flat. Consider the corresponding holomorphic bundle $F_\alpha$
on $\Tw(T)$. Clearly, the $(0,1)$-part of $\nabla_\alpha$,
taken with respect to $L$, is equal to $\bar\6_{gr}+\nu$.
Therefore, $F_\alpha\restrict{(T,L)} \cong B$,
for all $\alpha$. To prove \ref{_C^tw_I_isomo_Theorem_},
it remains to show that $F_\alpha \in \c C^{tw}_I$ for 
some $\alpha$. Let $(\cdot)^{0,1}_I$ denote the
operation of taking the $(0,1)$-Hodge component
with respect to $I$. Clearly, 
\begin{equation}\label{_C^tw_I_in_terms_of_nu_alpha_Equation_}
F_\alpha \in \c C^{tw}_I \iff  (\nu_\alpha)^{0,1}_I=0.
\end{equation}
An elementary calculation insures that
if $I\neq \pm L$, then
\[
x (\nu)^{0,1}_I = y (\bar\nu)^{0,1}_I,
\]
where $x, y$ are non-zero complex numbers,
depending on $I, L$. Choosing $\alpha := - \frac x y $, we
obtain that $(\nu_\alpha)^{0,1}_I=0$. By 
\eqref{_C^tw_I_in_terms_of_nu_alpha_Equation_},
this means that $F_\alpha \in \c C^{tw}_I$. We have
shown that the full and faithful functor
\[ R_L:\; \c C^{tw}_I \arrow \c C^b(T,L)
\]
is surjective on the set of isomorphism classes of objects.
Therefore, $R_L$ is an equivalence. We proved \ref{_C^tw_I_isomo_Theorem_}
for $M$ a generic torus of complex dimension $ >2$.

\subsection{Reflexive sheaves on hyperk\"ahler surfaces}
\label{_R_L_equi_refle_for_surfa_Subsection_}

In this Subsection, we prove \ref{_C^tw_I_isomo_Theorem_}
for $M$ a 2-dimensional torus or a K3 surface. 
Fix a hyperk\"ahler structure on $M$,
and let $L$ be a generic induced complex structure.
We have identified the category $C^r(M,L)$
with the direct limit \[ \lim_\arrow \c C(M\backslash S_i, L). \]
where $S_i$ is a sistem of finite subsets of $M$
ordered by inclusion, and $\c C(M \backslash S_i)$ the 
category of coherent sheaves on $M \backslash S_i$ 
which can be extended to coherent sheaves on $M$. 
If $S$ is non-empty, then the category of coherent sheaves
on $(M\backslash S_i, L)$ has cohomological dimension 
$\leq 1$ (\ref{_open_coho_dim_less_n_Corollary_}).
Therefore, $C^r(M,L)$ has cohomological dimension 
$\leq 1$ as well. Simple objects of this category
are stable bundles on $(M,L)$. The simple objects
of $\c C^{tw}_I$ are bundles
of form $\Tw(B)$ on $\Tw(M)$; therefore,
the simple objects of $\c C^{tw}_I$ 
in 1-to-1 correspondence with the
simple (that is, stable) reflexive sheaves
on $(M,L)$.

Consider the category 
\[ \c C^r:= \lim_\arrow \c C^b(\Tw(M \backslash S))
\]
(see Subsection \ref{_C^tw_I_defini_Subsection_}).
Let $\c C^{tw}\subset \c C^r$
be a full subcategory of $\c C^r$ consisting
of all extensions of coherent sheaves which 
have locally free reflexizations of form $\Tw(B)$,
where $B$ is a hyperholomorphic bundle. Clearly,
$\c C^{tw}_I$ is a full subcategory of $\c C^{tw}$. By 
\ref{_dire_on_tw_with_suppo_compu_Proposition_},
we have 
\[ \Ext^2_{\c C^r}(\Tw(B), \Tw(B'))=0
\]
for all hyperholomorphic bundles $B, B'$.
Therefore, in $\c C^{tw}$, the second Yoneda extension
between simple objects vanishes. By \ref{_fini_length_coho_1_Lemma_},
$\c C^{tw}$ has cohomoloical dimension $\leq 1$.
Now \ref{_cohomo_dime_subcate_Lemma_}
implies that $\c C^{tw}_I$ has 
cohomological dimension $\leq 1$ as well.

We arrive to assumptions
of \ref{_functo_inducing_isomo_on_semisi_Proposition_}:
a functor of abelian categories of finite length and
cohomological dimension $\leq 1$ induces an isomorphism
on the respective categories of semisimple objects.

To prove \ref{_C^tw_I_isomo_Theorem_} it remains to show
that $R_L:\c C^{tw}_I \arrow \c C^r(M,L)$ induces an isomorphism
on the first Yoneda extensions between simple objects.

Let $\nu$ be an Yoneda extension from $\c C^{tw}_I$,
\begin{equation}
\nu \in \Ext^1_{\c C^{tw}_I}(F_1, F_2),
\end{equation}
where $F_1, F_2\in \c C^{tw}_I$ are semisimple objects.
Then $F_i = \Tw(B_i)$, where $B_i$ is a hyperholomorphic
bundle on $M$. Furthermore, $\nu$ is an element of 
\[ H^0 \bigg(R^1\pi_* \bigg(
   F_1^* \otimes_{\calo_{\Tw(M)}} F_2\restrict{\Tw(M\backslash S)}
\bigg)\bigg)
\]
where $S$ is a finite set. The restriction
of $\nu$ to $(M\backslash S, I)\subset \Tw(M\backslash S)$
vanishes, 
and the image $r(\nu)$ of $\nu$ inside 
\begin{equation}\label{_kernel_determines_local_coho_Equation_}
\begin{aligned} 
P := \ker\bigg( 
\oplus_{x_i\in S} & \Hom_\C(\hat\calo_{s_{x_i}},\C)\otimes_{\calo_{\Tw(M)}}
  \Tw(B_1^* \otimes B_2) \otimes_{\calo(\C P^1)} \calo(2) \\
  {}& \arrow R^2\pi_* \Tw(B_1^* \otimes B_2) \bigg)
\end{aligned}
\end{equation}
is trivialized over $\C P^1\backslash I$ 
as in \eqref{_loca_coho_sheaf_trivialized_Equation_}
(see \ref{_exte_in_C^tw_I_in_terms_of_loca_Claim_}).

By \ref{_dire_on_tw_with_suppo_compu_Proposition_},
we have an exact sequence
\begin{align*}
0 \arrow & \Ext^1(B_1, B_2)\otimes H^0_{\C P^1}(\calo(1)) \arrow
  \Ext^1(F_1\restrict{\Tw(M\backslash S)},
  F_2\restrict{\Tw(M\backslash S)})  \\
 {}& \stackrel r\arrow P \arrow 0.
\end{align*}
An extension of reflexive  sheaves 
\[ \nu\in \Ext^1(F_1\restrict{\Tw(M\backslash S)},
   F_2\restrict{\Tw(M\backslash S)})
\]
belongs to $\c C^{tw}_I$ iff
$r(\nu)$ is trivialized over $\C P^1 \backslash I$
and $\nu \restrict{(M, I)}=0$.
Since $r(\nu)$ is trivialized over $\C P^1\backslash I$,
the restriction
$r(\nu)\restrict{\C P^1\backslash I}$ is determined by its value at any point
$L\in \C P^1$, $L\neq I$. Similarly, any section
of \[ \Ext^1(B_1, B_2)\otimes H^0_{\C P^1}(\calo(1))\]
vanishing at $I$ is determined by its restriction
on $(M,L)$. Therefore, $R_L:\; \c C^{tw}_I \arrow \c C^r(M,L)$
induces an embedding on $\Ext^1$ between semisimple objects.
To show that this map is an isomorphism, take
\[ \alpha \in \Ext^1_{\c C^b(M\backslash S, L)}(B_1, B_2), \]
where $B_1, B_2$ are hyperholomorphic bundles considered
as holomorphic bundles over $(M, L)$. 
By \ref{_B^star_trivialization_extended_Proposition_},
any trivialized over $\C P^1\backslash I$
section of $P\restrict{\C P^1\backslash I}$ 
can be extended to $\C P^1$. Therefore, the image
\begin{align*}
   r(\alpha)\in P_L& :=
   \ker\bigg( 
   \Hom_\C(\hat\calo_{(M,L), x},\C)\otimes_{\calo_{(M,L)}}K(M,L)
   \otimes B_1^* \otimes B_2\\{}&\arrow 
   H^2((M, L), B_1^* \otimes B_2)\bigg)
\end{align*}
 can be extended to a section of $P$, on the whole
$\C P^1$.  Clearly, any element in $\Ext^1(B_1, B_2)$
can be extended to a unique section of
$\Ext^1(B_1, B_2)\otimes H^0_{\C P^1}(\calo(1))$
vanishing at $I$. We have a commutative diagram of long exact sequences,
with the rightmost and the leftmost vertical arrows isomorphisms

{\scriptsize
\minCDarrowwidth1.5pc
\begin{equation}\label{_squa_long_exa_R_L_Equation_}
\begin{CD}
0 @>>> \Ext^1(B_1, B_2)\otimes 
       H^0_{\C P^1}(\calo(1))_{I}@>>> 
       \Ext^1_{\c C^{tw}_I}(F_1, F_2) @>>> P_{tr}  @>>> 0\\
&& @VVV @V {R_L} VV @VVV&&\\
0 @>>> \Ext^1(B_1, B_2) @>>> 
       \Ext^1_{\c C^b(M\backslash S, L)}(B_1, B_2) 
       @>>> P_L   @>>> 0
\end{CD}
\end{equation}
}
where $H^0_{\C P^1}(\calo(1))_{I}$ denotes the space
of sections of  $\calo(1)$ vanishing at $I$, 
$P_{tr}$ denotes the space of sections of
$P$ which are trivialized over $\C P^1\backslash I$,
and the vertical maps are restrictions to
$(M,L) \subset\Tw(M)$.

The middle vertical arrow of \eqref{_squa_long_exa_R_L_Equation_}
is also an isomorphism, by the snake lemma;
therefore, $R_L$ induces an isomorphism on 
Yoneda extensions of semisimple objects. This proves
\ref{_C^tw_I_isomo_Theorem_}
for $M$ a 2-dimensional torus or a K3 surface.
\ref{_C^tw_I_isomo_Theorem_} is proven.
\endproof


\section{Coherent sheaves with isolated singularities}
\label{_cohe_isola_Section_}


The main result of this Section is the following

\hfill

\theorem\label{_c_C_L_independent_Theorem_}
Let $M$ be a hyperk\"ahler K3 surface or a 
compact torus without 
non-trivial integer $SU(2)$-invariant classes, and $L$, $L'$
generic induced complex structures. Consider the categories $\c C(M,L)$,
$\c C(M,L')$ of coherent sheaves on $(M, L)$ and $(M,L')$. 
Then the categeries $\c C(M,L)$,
$\c C(M,L')$ are equivalent. Moreover, this equivalence
can be chosen canonically if we fix an induced complex
structure $I \neq L, L'$.

\hfill

{\bf Proof:} 
Let $F\in \c C(M,L)$, and $F^{**}$ its reflexization, 
which is a  bundle. 
Fix an induced complex
structure $I \neq L, L'$.
As in the proof of \ref{_refle_shea_isomo_Theorem_},
we consider the category 
$\c C^{tw}_I$. We have shown that
the restriction $R_L:\; \c C^{tw}_I\arrow \c C^r(M,L)$
defines an equivalence, where $\c C ^r(M,L)$
is the category of reflexive sheaves
on $(M,L)$ (\ref{_C^tw_I_isomo_Theorem_}).
Denote by $E\in \c C^{tw}_I$
the reflexive sheaf which corresponds to
$F^{**}$ under this equivalence.

Our argument uses the local twistor 
geometry of a hyperk\"ahler manifold, following
Subsection \ref{_twi_spa_local_Subsection_}
(see also \cite{_Verbitsky:hypercomple_}).
The following claim is trivial.

\hfill

\claim\label{_C^tw_I_trivialized_Claim_}
In the above assumptions, let
$E\in \c C^{tw}_I$, and let $x$ be any point in $M$.
Then the bundle $E$ is canonically 
 trivialized in a neighbourhood
of $s_x \backslash I$, where 
$s_x \subset \Tw(M)$ is a horisontal
twistor section.

\hfill

{\bf Proof:}
Let $x$ be a point in $A$, where $A\times \{ I\}$ is a singular set
of $E$ which is used in the definition of
$\c C^{tw}_I$. 
In a neighbourhood of $s_x \backslash I$, the sheaf
$E$ is trivialized, by definition of $\c C^{tw}_I$. 
If $x\notin A$, then $E$ in a neighbourhood
of $s_x$ is a successive extension
of the bundles of form $\Tw(B_i)$.
Each of these bundles is equal to
$\calo_{\C P^1}^{\oplus k_i}$ on $s_x$,
and therefore, $E\restrict{s_x}$ is also
trivial. Using \eqref{_E_loca_triv_Equation_}, we 
obtain a trivialization of $E$ in a neighbourhood 
of $s_x \backslash I$, for all $x\in M$. 
\endproof

\hfill

We are going to attach singularities to $E$ in such a way that
the equivalence $R_L:\; \c C^{tw}_I\arrow \c C ^b(M,L)$
can be extended to the sheaves with isolated
singularities.

\hfill

We prove \ref{_c_C_L_independent_Theorem_}
in the same fashion as we proved \ref{_refle_shea_isomo_Theorem_}.
Just as in \ref{_C^tw_I_isomo_Theorem_}, we shall
construct an intermediate category $\tilde{\c C}^{tw}_I$,
which is a subcategory of the category of 
coherent sheaves 
on the twistor space. We show that the 
natural restriction functor from $\tilde{\c C}^{tw}_I$
to $\c C(M,L)$ is an equivalence, if $L$ is a generic
complex structure distinct from $I$. 
Since $L$ can be chosen arbitrarily,
this will imply immediately that
$\c C(M,L)$ is equivalent to $\c C(M,L')$.

\hfill

\definition\label{_tilde_C(I)_Definition_}
The category $\tilde{\c C}^{tw}_I$
is defined as follows. 

An object of $\tilde{\c C}^{tw}_I$
is a coherent sheaf $F$ on $\Tw(M)$
satisfying the following conditions.

\begin{description}
\item[(i)] The reflexization $E:= F^{**}$
belongs to $\c C^{tw}_I$,
where the category $\c C^{tw}_I$ is defined 
as in \ref{_C_I^tw_Definition_}.

\item[(ii)] The sheaf $F$ is non-singular outside
of $A\times \C P^1\subset \Tw(M)$,
where $A\subset M$ is a finite set.

\item[(iii)] Let $Z:= \{A\} \times  I \subset (M, I) \subset \Tw(M)$
be the finite subset of $(M, I) \subset \Tw(M)$
corresponding to $A$, and 
$j:\; \Tw(M)\backslash Z \hookrightarrow \Tw(M)$
the natural embedding. Then the canonical
homomorphism $F\arrow j_* j^* F$ is an
isomorphism.

\item[(iv)] In the notation of of Subsection
\ref{_twi_spa_local_Subsection_}, 
consider the restriction of $F$ and of $E= F^{**}$
to $\Tw(M)_{x_0, I}$. Let
$E\restrict {\Tw(M)_{x_0, I}}\stackrel {\tau_E}\arrow  \xi^*E_S$
be the trivialization constructed in
\ref{_C^tw_I_trivialized_Claim_}.
Then there exists a sheaf $F_S$ 
equipped with a map $F_S\stackrel{\rho_{F_S}}\arrow E_S$
which is equivalence outside of $x_0$,\footnote{This means that $\rho_{F_S}$ is
the reflexization map.}
and an isomorphism 
\[ F\restrict {\Tw(M)_{x_0, I}}\stackrel{\tau_F}\arrow  \xi^*F_S
\]
such that the following diagram is commutative
\[
\begin{CD}
F\restrict {\Tw(M)_{x_0, I}} @>{\tau_F}>> \xi^*F_S\\[3mm]
@V{\rho_F}VV @VV{\xi^* (\rho_{F_S})}V \\[3mm]
E\restrict {\Tw(M)_{x_0, I}} @>{\tau_E}>> \xi^*E_S
\end{CD}
\]
(Here, $\rho_F$ and $\rho_{F_S}$ denote
the respective reflexization maps).
This means that the singularities of 
$F$ are compatible with the local
trivialization of $E$
constructed in Subsection
\ref{_twi_spa_local_Subsection_}.
\end{description}

The morphisms of $\tilde{\c C}^{tw}_I$
are morphisms of coherent sheaves.

\hfill

\remark
Clearly, the sheaf $F_S$ and the morphism
$F_S\arrow E_S$ can be reconstructed in a canonical way
from $F$. Take for $S$ the neighbourhood of $x_0$ in $(M, L)$.
Then $E_S$ is the restriction of
$E$ to $S\subset (M, L)\subset \Tw(M)$,
and $F_S\cong  F\restrict S$. 

\hfill

The following theorem implies
\ref{_c_C_L_independent_Theorem_},
as we indicated in the beginning of this Section.

\hfill

\theorem\label{_tilde_C(I)_C_I_equiv_Theorem_}
Let $M$ be a compact hyperk\"ahler manifold, $I$ an arbitrary induced
complex structure, and $\tilde{\c C}^{tw}_I$ the category
constructed above. Consider a generic induced complex structure
$L$, and let 
\begin{equation}\label{_restri_tilde_C(I)_Equation_}
R_L:\; \tilde{\c C}^{tw}_I\arrow \c C(M,L)
\end{equation}
be the restriction map, $F \arrow F\restrict{(M,L)}$.
Then \eqref{_restri_tilde_C(I)_Equation_}
is an equivalence of categories.

\hfill

{\bf Proof:} We construct the inverse functor 
$\Upsilon:\;\c C(M,L)\arrow  \tilde{\c C}^{tw}_I$
as follows. Take $F_L\in \c C(M,L)$.
Let $E_L$ be its reflexization,
and $E\in \c C^{tw}_I$ the corresponding
sheaf over the twistor space (\ref{_C^tw_I_isomo_Theorem_}). 
Let $x_0\in M$ be a singular point of $F_L$.
Consider the infinitesimal neighbourhood
of 
\[ \{x_0\}\times(\C P^1\backslash I)\subset \Tw(M) \backslash (M,I)
\] in $\Tw(M) \backslash (M,I)$, denoted, as in Subsection
\ref{_twi_spa_local_Subsection_},
 by $\Tw(M)_{x_0, I}$.
Then $E$ is trivialized over $\Tw(M)_{x_0, I}$,
and we may write $E\restrict{\Tw(M)_{x_0, I}} = \xi^*(E_{L, x_0})$,
where $E_{L, x_0}$ is the germ of $E_L$ in
$x_0$. Let $F\restrict{\Tw(M)_{x_0, I}}:=\xi^*(F_{L, x_0})$ be 
the sheaf corresponding to $F_L$.
Since $F_L$ has isolated singularities, the sheaves
$F\restrict{\Tw(M)_{x_0, I}}$ and
$E\restrict{\Tw(M)_{x_0, I}}$ are 
canonically isomorphic outside of
$\{x_0\}\times(\C P^1\backslash I)$.
Gluing together $F\restrict{\Tw(M)_{x_0, I}}$
and $E$, 
we obtain a sheaf $F_0$ on $\Tw(M) \backslash (M,I)$.
Let $A$ be the singular set of $F_L$. Then
$E$ is equal to $F_0$ outside of 
$\{A\} \times \C P^1\subset \Tw(M)$,
hence $F_0$ can be extended smoothly
from $\Tw(M) \backslash (M,I)$
to a sheaf $F_1$ on $\Tw(M) \backslash Z$, where
$Z= \{A\} \times \{I\}$ is the finite set
defined in \ref{_tilde_C(I)_Definition_}.
Let $\Upsilon(F_L) := j_*(F_1)$, where 
$j:\; \Tw(M) \backslash Z \hookrightarrow \Tw(M)$
is the natural embedding.
Since $Z$ has codimension $\geq 3$,
the direct image $\Upsilon(F_L)= j_*(F_1)$
is a coherent sheaf, and we have
$R_L (\Upsilon(F_L)) = F_L$. By construction,
this sheaf satisfies the conditions
of \ref{_tilde_C(I)_Definition_}.
This construction is clearly functorial.
Indeed, it is functorial on bundles
by \ref{_C^tw_I_isomo_Theorem_},
and the sheaf $F_L$ is just a bundle
with some singularities attached;
every morphism of sheaves from $\c C_I$ induces
a morphism on the corresponding reflexization
bundles, which is compatible with the
trivialization \eqref{_E_loca_triv_Equation_},
hence extends naturally to
isolated singularities.

We proved \ref{_tilde_C(I)_C_I_equiv_Theorem_}.
\ref{_c_C_L_independent_Theorem_} is proven.
\endproof

\hfill

From \ref{_generic_are_connected_Theorem_}
and \ref{_c_C_L_independent_Theorem_}, the following result
is apparent.

\hfill

\theorem\label{_generic_are_equi_Theorem_}
Let $M$ be a compact torus or a K3 surface, 
and $L$, $L'$ complex structures of K\"ahler type. 
Assume that $L$ and $L'$ are Mumford-Tate generic. 
Then the categories $\c C(M,L)$, $\c C(M,L')$
of coherent sheaves on $(M, L)$ and $(M,L')$
are equivalent.

\hfill

{\bf Proof:} Consider the sequence $\c H_1, \c H_2, ... \c H_n$
of hyperk\"ahler structures, satisfying the
assumptions of \ref{_generic_are_connected_Theorem_},
and let $L_0 = L, L_2, L_3, ..., L_{n+1} = L'$
be the corresponding sequence of Mumford-Tate generic complex
structures. By \ref{_MT_gene_gene_wrt_hk_Claim_},
the complex structures $L_k$ are generic.
Applying 
\ref{_c_C_L_independent_Theorem_}, we find that
the following categories are equivalent:
\[  \c C(M, L_k) \sim \c C(M, L_{k+1})
\]
Therefore, $\c C(M,{L_0})=\c C(M,L)$
is equivalent to $\c C(M, {L_{n+1}}= \c C(M,L')$.
We proved \ref{_generic_are_equi_Theorem_}. \endproof

\hfill

\hfill

{\bf Acknowledgements:} I am grateful to V. Ginburg, who
asked about the categories of coherent sheaves on a K3 surface, 
and D. Kaledin for many hours of stimulating discussions. Many thanks 
to I. Dolgachev, D. Kazhdan, A. Tyurin and A. Vaintrob for 
insightful comments and useful advice.

\hfill


\end{document}